\documentclass[12pt] {article}
\usepackage{amsmath,amsthm}
\usepackage{amssymb,latexsym}
\usepackage{graphicx}
\usepackage{amssymb}
\usepackage{amsfonts}
\usepackage{amscd}

\usepackage{amscd}
\title{Theory of valuations on manifolds, III. Multiplicative structure in the general case.}
\date{}
\author{ Semyon Alesker\footnote{Partially supported by ISF grant 1369/04.}, Joseph H.G. Fu \footnote{Partially supported by NSF grant DMS-0204826.}}

\def\eps{\varepsilon}
\def\alp{\alpha}
\def\ome{\omega}
\def\Ome{\Omega}
\def\lam{\lambda}

\def\to{\rightarrow}
\def\qed { Q.E.D. }
\def\pt{\partial}

\def\RR{\mathbb{R}}
\def\CC{\mathbb{C}}

\def\NN{\mathbb{N}}
\def\ZZ{\mathbb{Z}}

\def\PP{\mathbb{P}}

\swapnumbers
\newtheorem{theorem}{Theorem}[subsection]
\newtheorem{corollary}[theorem]{Corollary}
\newtheorem{lemma}[theorem]{Lemma}
\newtheorem{proposition}[theorem]{Proposition}

\theoremstyle{definition}

\newtheorem{definition}[theorem]{Definition}

\theoremstyle{proposition-definition}
\newtheorem{proposition-definition}[theorem]{Proposition-Definition}

\def\ca{{\cal A}} \def\cb{{\cal B}} 
\def\cd{{\cal D}}  
 \def\ch{{\cal H}} 
 \def\ck{{\cal K}} 
  
\def\cp{{\cal P}} 

 \def\ct{{\cal T}} 
\def\cv{{\cal V}}  
 
\def\mov{|\omega_V|}
\def\inj{\hookrightarrow}
\def\surj{\twoheadrightarrow}
\def\Nor{\operatorname{Nor}}
\def\reach{\operatorname{reach}}

\newcommand\mass{\operatorname{mass}}
\newcommand \supp{\operatorname{supp} \,}
\newcommand \dist{\operatorname{dist} \,}

\begin{document}
\maketitle \setcounter{section}{-1}
\begin{abstract}
This is the third part of a series of articles where the theory of
valuations on manifolds is constructed. In \cite{part2} the notion
of a smooth valuation on a manifold was introduced. The goal of
this article is to put a canonical multiplicative structure on the
space of smooth valuations on general manifolds, thus extending
some of the affine constructions from \cite{alesker-poly},
\cite{part1}.
\end{abstract}
\tableofcontents
\section{Introduction.}
In convexity there are many geometrically interesting and well
known examples of valuations on convex sets, including Lebesgue measure, the
Euler characteristic,  surface area, mixed volumes, and  affine
surface area. For a description of older classical developments on
this subject we refer to the surveys \cite{mcmullen-schneider},
\cite{mcmullen-survey}. For  general background on convexity we
refer to the book \cite{schneider-book}.

Approximately during the last decade there has been significant
progress in this classical subject, leading to new
classification results for various classes of valuations, and to the
discovery of new structures on them. This progress has shed new
light on the notion of valuation, permitting an extension to the
more general setting of valuations on manifolds, and to sets  which are not necessarily convex (a concept which in any case has no meaning on a general
manifold). The development of the theory of valuations on
manifolds was started by one of the authors in the first two parts  \cite{part1}, \cite{part2}
of the present series of articles.

In \cite{part2} the notion of smooth valuation on a smooth
manifold was introduced. Roughly put a smooth valuation can be
thought as a finitely additive $\CC$-valued measure on a class of
nice subsets; this measure is required to satisfy some additional
assumptions of continuity (or rather smoothness) in some sense.
The basic examples of smooth valuations on a general manifold $X$
are smooth measures on $X$ and the Euler characteristic. Moreover
the well known intrinsic volumes of convex sets can be generalized
to provide examples of smooth valuations on an arbitrary {\itshape
Riemannian} manifold; these valuations are known as
Lipschitz-Killing curvatures.

The goal of this article is to put a canonical multiplicative
structure on the space of smooth valuations on a general smooth manifold.
When the manifold is an affine space the multiplicative structure
was defined in \cite{part1} (in the more specific situation of
valuations polynomial with respect to translations it was defined
in \cite{alesker-poly}). The construction of the product on
general manifolds presented in this article uses the construction
from \cite{part1} for the affine case in combination with
tools from geometric measure theory. Roughly it works as
follows. Choosing a coordinate atlas for $X$, one uses the product of valuations on
$\RR^n$, defined by the construction of \cite{part1}, to define the product locally. Then one
shows that the products obtained on each coordinate patch
coincide on pairwise intersections, and that the result does not
depend on the choice of atlas. This step uses geometric measure theory.


Let us denote by $V^\infty(X)$ the space of smooth valuations on
$X$. The product
$$V^\infty(X)\times V^\infty(X)\to V^\infty(X)$$
defined in this article is a continuous map, with respect to which $V^\infty(X)$ becomes a commutative associative
algebra with unit (where the unit is the Euler characteristic).

In \cite{part2} it was shown that the assignment to any open subset $U\subset X$
$$U\mapsto V^\infty(U),$$
with the natural operations of restriction, defines a sheaf of
vector spaces on $X$ denoted by $\cv^\infty_X$. The multiplicative
structure on smooth valuations defined in this article commutes
with restriction to open subsets. Hence $\cv^\infty_X$ becomes a
sheaf of commutative associative algebras with unit. Further
properties of the multiplicative structure are studied in the
fourth part of this series \cite{part4}.

The article is organized as follows. Section \ref{background}
contains some background; it does not contain new results. There
we also fix some notation used throughout the article. In Section
\ref{normalcycles} we discuss normal cycles; the exposition
follows mostly \cite{Fu:1989}-\cite{fu:1994}. Then we explain how
to use normal cycles to construct valuations. We also recall
related results from geometric measure theory. In Section
\ref{auxiliary} we prove some auxiliary results of  a technical
nature. The main part of the article is Section
\ref{construction}, where we present a construction of the product
on smooth valuations and prove that it is independent of the
choices involved.

{\bf Acknowledgements.} The first named author is grateful to J. Bernstein for
very useful discussions, and to V.D. Milman for his interest in this work and useful
discussions.

\section{Background.}\label{background}

In Subsection \ref{notation} we fix some notation which will be
used throughout the article. In Subsection \ref{convexity} we
recall some results from \cite{alesker-int} and deduce from them
Corollary \ref{cor-mink}, which will be used later on. In Subsection
\ref{valuation_theory} we recall some necessary facts from the theory
of valuations and some results from \cite{part1} and \cite{part2}.
\subsection{Notation.}\label{notation}
Let $V$ be a finite dimensional real vector space.

$\bullet$ Let $\ck(V)$ denote the family of convex compact subsets of $V$.

$\bullet$ Let $\RR_{\geq 0}$ denote the the set of non-negative real numbers.

$\bullet$ For a manifold $X$ let us denote by $|\ome_X|$ (or just by $|\ome|$ if it does not lead to confusion)
the line bundle of densities over $X$.

$\bullet$ For a smooth manifold $X$ let $\cp(X)$ denote the family of all simple subpolyhedra of $X$.
(Namely $P\in \cp(X)$ iff $P$ is a compact subset of $X$ locally diffeomorphic to
$\RR^k\times \RR^{n-k}_{\geq 0}$ for some $0\leq k\leq n$.
For a precise definition see \cite{part2}, Subsection 2.1.)

$\bullet$ We denote by $\PP_+(V)$ the {\itshape oriented
projectivization} of $V$. Namely $\PP_+(V)$ is the manifold of
oriented lines in $V$ passing through the origin.

$\bullet$ For a convex compact set $A\in \ck(V)$ let us denote by
$h_A$ the {\itshape supporting functional} of $A$, $h_A\colon
V^*\to \RR$. It is defined by
$$h_A(y):=\sup\{y(x)| x\in A\}.$$

$\bullet$ Let $L$ denote the (real) line bundle over $\PP_+(V^*)$
such that its fiber over an oriented line $l\in \PP_+(V^*)$ is
equal to the dual line $l^*$.

$\bullet$ A subset $A$ of a Euclidean space $V$ is said to have
{\itshape positive reach} (or to be {\itshape semi-convex}) if there
exists $\eps >0$ such that for any $x\in V$ with $\dist(x,A)<\eps$
there exists a unique point $y\in A$ with $\dist(x,y)=\dist(x,A)$ (cf. \cite{Federer:1959}). The reach of $A$ is defined to be the supremum of all such $\eps$. Note that $A$ is convex iff $\reach(A) = \infty$.

\subsection{Some convexity.}\label{convexity}
First let us recall some results from \cite{alesker-int}.  Let
$\bar K=(K_1 , K_2 ,\dots, K_s )$ be an $s$-tuple of compact
convex subsets of $V $. Let $r\in \NN\cup \{\infty\}$. For any
$\mu \in C^{r}(V,\mov)$ consider the function $M_{\bar K} F \,
:\RR^{s}_{+} \to \CC \, ,\mbox{where }\RR_{+}^{s}=
\{(\lam_1,\dots,\lam_s) \,| \,\lam_j\geq0 \} $, defined by    $$ (
M_{\bar K} \mu )
  (\lam _1,\dots,\lam _s)= \mu(\sum_{i=1}^{s} \lam _i K_i) .  $$
\begin{theorem}[\cite{alesker-int}]\label{minkowski}
(1) $M_{\bar K}\mu \in C^{r}(\RR^s_+)$ and $M_{\bar K}$ is a
continuous operator from $C^{r}(V,\mov)$ to $C^{r}(\RR^s_+)$.

(2) Assume that  a sequence $\mu^{(m)}$ converges to $\mu$ in
$C^{r}(V,\mov)$. Let $ K_{i}^{(m)}, \, K_{i},\,
i=1,\dots,s,\,m\in\NN $ be convex compact sets in $ V $, and for
every $ i=1,\dots ,s $ $K_{i}^{(m)} \to K_{i} $ in the Hausdorff
metric as $m\to \infty$. Then $ M_{\bar K^{(m)} } \mu^{(m)} \to
M_{\bar K } \mu $ in $C^{r}(\RR^s_+)$ as $m\to \infty$.
\end{theorem}

\begin{corollary}\label{cor-mink}
Let $r\in \NN$. Let $R>0$. Then there exists a constant $C$
depending on $r,R,$ and $n$ only such that for any $\mu\in
C^r(V,|\ome_V|)$ and any $K\in \ck(V)$ contained in a centered
Euclidean ball of radius $R$ one has
\begin{eqnarray*}
\big| \frac{\pt^r}{\pt\lam_1\dots\pt\lam_r}\big|_0
\mu(K+\sum_{i=1}^r\lam_iA_i)\big|\leq C\cdot ||\mu||_{C^r(C R\cdot
D)} \cdot \prod_{i=1}^r ||h_{A_i}||_{C^2(\PP_+(V^*))}.
\end{eqnarray*}
\end{corollary}
{\bf Proof.} Consider the functional
\begin{eqnarray*}
\psi\colon C^r(V,|\ome_V|)\times \ck(V)^{r+1}\to \CC
\end{eqnarray*}
given by
\begin{eqnarray*}
\psi(\mu;K;A_1,\dots,A_r)=\frac{\pt^r}{\pt\lam_1\dots\pt\lam_r}\big|_0
\mu(K+\sum_{i=1}^r\lam_iA_i).
\end{eqnarray*}
By Theorem \ref{minkowski} $\psi$ is a continuous map. Clearly $\psi$ is linear with respect
to $\mu$ and symmetric with respect to $A_1,\dots,A_r$. Moreover $\psi$ is Minkowski additive
with respect to each of $A_1,\dots,A_r$. (Minkowski additivity, say with respect to $A_1$, means
that
\begin{eqnarray*}
\psi(\mu;K;\alp A_1'+\beta A_1'',A_2,\dots,A_r)=
\alp \psi(\mu;K; A_1',A_2,\dots,A_r)+\beta \psi(\mu;K; A_1'',A_2,\dots,A_r)
\end{eqnarray*}
for any $\alp,\beta\geq 0,\, A_1',A_1''\in \ck(V)$.)

Let $h\in C^2(\PP_+(V^*),L)$. Then $h$ can be presented as
\begin{eqnarray}\label{z1}
h=h_{A'}-h_{A''}
\end{eqnarray}
where $h_{A'},h_{A''}\in C^2(\PP_+(V^*),L)$ are supporting functionals
of convex compact sets $A',A''\in \ck(V)$ respectively, and
\begin{eqnarray}\label{z2}
\max\{ ||h_{A'}||_{C^2(\PP_+(V^*))},||h_{A''}||_{C^2(\PP_+(V^*))}\}\leq
c ||h||_{C^2(\PP_+(V^*))}
\end{eqnarray}
where $c$ is a constant depending on $n$ only. Indeed let us take
$h_{A'}=h+T\cdot h_D,\, A''=T\cdot D$ where $D$ is the unit Euclidean ball,
and $T>0$ a large enough constant depending on $||h||_{C^2(\PP_+(V^*))}$.
Now let us extend $\psi$ to a functional
\begin{eqnarray*}
\tilde\psi\colon C^r(V,|\ome_V|)\times \ck(V)\times
(C^2(\PP_+(V^*),L))^r\to\CC
\end{eqnarray*}
by linearity. More precisely let $(\mu;K;h_1,\dots,h_r)\in
C^r(V,|\ome_V|)\times\ck(V) \times (C^2(\PP_+(V^*),L))^r$. Let us
define $\tilde\psi(\mu;K;h_1,\dots,h_r)$ recursively on the number
of non-convex functions among $h_1,\dots,h_r$. If this number is
equal to zero, there is nothing to define. Assume we have defined
$\tilde\psi$ for $k-1<r$ not necessarily convex functions. Let us
define it for $k$ such functions. We may assume that $h_i,\, i>k$ are
convex. Choose a presentation
\begin{eqnarray}\label{z3}
h_k=h_{A_k'}-h_{A_k''}
\end{eqnarray}
as in (\ref{z1}), and satisfying (\ref{z2}). Now define
\begin{eqnarray*}
\tilde\psi(\mu;K;h_1,\dots,h_{k-1},h_k,h_{k+1},\dots,h_r)=\\
\tilde\psi(\mu;K;h_1,\dots,h_{k-1},h_{A_k'},h_{k+1},\dots,h_r)-
\tilde\psi(\mu;K;h_1,\dots,h_{k-1},h_{A_k''},h_{k+1},\dots,h_r)
\end{eqnarray*}
where the right hand side is defined by the assumption of induction.

It is easy to see that the extension $\tilde \psi$ is well defined (i.e.
it does not depend on the choice of presentation (\ref{z3})). Next $\tilde\psi$
is continuous due to (\ref{z2}), and it is linear with respect to $h_1,\dots,h_r$.
Now Corollary \ref{cor-mink} follows from a very general, simple, and well known lemma
as follows.
\begin{lemma}\label{known-lemma}
Let $X$ be a compact topological space. Let $F_1,\dots,F_t$ be locally convex
$\CC$-linear topological spaces. Let
$$\phi\colon X\times F_1\times \dots \times F_t\to \CC$$
be a continuous map which is linear with respect to the last $t$ arguments.

Then there exist continuous semi-norms $||\cdot||_1,\dots,||\cdot||_t$ on
$F_1,\dots,F_t$ respectively such that for any $x\in X, \xi_1\in F_1,\dots,\xi_t\in F_t$
one has
\begin{eqnarray*}
|\phi(x,\xi_1,\dots,\xi_t)|\leq \prod_{i=1}^t||\xi_i||_i.
\end{eqnarray*}
\end{lemma}
Thus Corollary \ref{cor-mink} is proved. \qed

\subsection{Some valuation theory.}\label{valuation_theory}
\begin{definition}
a) A function $\phi :{\cal K}(V) \to \CC$ is called a valuation if
for any $K_1, \, K_2 \in {\cal K}(V)$ such that their union is
also convex one has
$$\phi(K_1 \cup K_2)= \phi(K_1) +\phi(K_2) -\phi(K_1 \cap K_2).$$

b) A valuation $\phi$ is called continuous if it is continuous
with respect to the Hausdorff metric on ${\cal K}(V)$.
\end{definition}

For the classical theory of valuations we refer to the surveys
McMullen-Schneider \cite{mcmullen-schneider} and McMullen
\cite{mcmullen-survey}. For general background from convexity
we refer to Schneider \cite{schneider-book}.

In \cite{part1} one has introduced a class $SV(V)$ of valuations called
{\itshape smooth valuations}. We refer to \cite{part1} for an axiomatic definition. Here
we only mention that $SV(V)$ is a $\CC$-linear space (with the obvious operations) with a
natural Fr\'echet topology. In this article we will need another description of $SV(V)$, given in Theorem \ref{onto}
below.

Let us denote by ${}\!^ {\textbf{C}}L$ the (complex) line bundle
over $\PP_+(V^*)$ whose fiber over $l\in \PP_+(V^*)$ is equal to
$l^*\otimes_{\RR}\CC$ (where $l^*$ denotes the dual space to $l$).

Note that for any convex compact set $A\in \ck(V)$ the supporting
functional $h_A$ is a continuous section of  ${}\!^ {\textbf{C}}L$:
$h_A\in C(\PP_+(V^*),{}\!^ {\textbf{C}}L)$.
\begin{theorem}[\cite{part1}, Corollary 3.1.7]\label{onto}
There exists a continuous linear map
\begin{eqnarray*}
\ct\colon \oplus_{k=0}^nC^\infty(V\times
\PP_+(V^*)^k,|\ome_V|\boxtimes {}\!^ {\textbf{C}}L^{\boxtimes k})
\to SV(V)
\end{eqnarray*}
which is uniquely characterised by the following property: for any
$k=0,1,\dots,n$, any $\mu\in C^\infty(V,|\ome_V|)$, any strictly
convex compact sets $A_1,\dots,A_k$ with smooth boundaries, and
any $K\in \ck(V)$ one has
\begin{eqnarray*}
\ct(\mu\boxtimes h_{A_1}\boxtimes \dots\boxtimes h_{A_k})(K)=
\frac{\pt^k}{\pt\lam_1\dots \pt\lam_k}\big|_0 \mu(K+\sum_{i=1}^k\lam_iA_i)
\end{eqnarray*}
where $\lam_i\geq 0$ in the right hand side.

Moreover the map $\ct$ is an epimorphism.
\end{theorem}

In \cite{part2} one has introduced for any smooth manifold $X$ a
class of finitely additive measures on the family of simple
subpolyhedra $\cp(X)$. This class is denoted by $V^\infty(X)$. It
is a $\CC$-linear space (with the obvious operations). Then
$V^\infty(X)$ has a natural Fr\'echet topology. Moreover in the
case  of a linear space $V$, the restriction of any element $\phi\in
V^\infty(V)$ to $\ck(V)\cap \cp(V)$ has a
(unique) extension by continuity in the Hausdorff metric to
$\ck(V)$, and this extension belongs to $SV(V)$. Thus one gets a
linear map
$$V^\infty(V)\to SV(V).$$
In \cite{part2}, Proposition 2.4.10, the following result was proved.
\begin{proposition}\label{isomorphism}
The  map $V^\infty(V)\to SV(V)$ constructed above is an isomorphism of
Fr\'echet spaces.
\end{proposition}



\newcommand{\lem}[1]{Lemma~\ref{lem:#1}}
\newcommand{\thm}[1]{Theorem~\ref{thm:#1}}
\newcommand{\prop}[1]{Proposition~\ref{prop:#1}}
\newcommand{\cor}[1]{Corollary~\ref{cor:#1}}

\newcommand\deriv{\frac{d}{dt}|_{t=0}}
\newcommand\G{\Gamma}
\newcommand\D{\Delta}
\newcommand\s{\sigma}
\newcommand\g{\gamma}
\newcommand\osc{\operatorname{osc}}
\newcommand\vol{\operatorname{vol}}
\newcommand\PD{\operatorname{PD}}
\newcommand\R{\mathbb R}
\newcommand\Rn{\mathbb R^n}
\newcommand\C{\mathbb C}
\newcommand\CP{\mathbb C \mathbb P}
\newcommand\Val{\operatorname {Val}}
\newcommand\ValSO{\operatorname {Val}^{SO(n)}}
\newcommand\ValU{\operatorname {Val}^{U(n)}}
\newcommand\Diff{\operatorname{Diff}}
\newcommand\nor{\operatorname{nor}}
\newcommand\tcone{\operatorname{tan}}
\newcommand\Tan{\operatorname{Tan}}
\newcommand\Hom{\operatorname{Hom}}
\newcommand\Vect{\operatorname{Vec}}
\newcommand\length{\operatorname{length}}
\newcommand\area{\operatorname{area}}
\newcommand\graph{\operatorname{graph}}
\newcommand\dvol{\operatorname{d\,vol}}
\newcommand\clos{\operatorname{clos}}
\newcommand\Endpt{\operatorname{End}}
\newcommand \ot{\otimes}
\newcommand \annih{\mathcal A}
\newcommand \K{\mathcal K}
\newcommand \lcur{[\![}
\newcommand \rcur{]\!]}

\section{Valuations and normal cycles.}\label{normalcycles}
\subsection{Normal cycles for convex sets.}
Let $V$ be a  real vector space with finite dimension $n$.
Let $K\in \ck(V)$. Let $x\in K$.
\begin{definition}\label{cc-1-1}
The tangent cone to $K$ at $x$ is the closure of the set $\{y\in V|\exists \eps>0\, \,
x+\eps y\in K\}$. We denote it by $T_xK$.
\end{definition}
It is easy to see that $T_xK$ is a closed convex cone.
\begin{definition}\label{cc-1-2}
The conormal cone to $K$ at $x$ is the set
$$\Nor^*(K,x):=\{y\in V^*| \,\, y(x)\leq 0 \forall x\in T_xK\}.$$
\end{definition}
Thus $\Nor^*(K,x)$ is also a closed convex cone. Put also
$$\Nor^*(A) :=  \cup_{x\in K}\left(\{x\}\times\Nor^*(K,x)\right).$$
Fixing a euclidean
metric $(, )$ on $V$, it will be convenient to define
${\Nor}(A)\subset V \times V$ as the image of
$\Nor^*(A)$ under the induced identification $V \times V^* \simeq
V \times V$. Finally, we put
$$\Nor_1(K): = \{(x,v)\in \Nor(K): |v| = 1\}. $$

It is easy to see that $\Nor^*(K)$ (resp. $\Nor(K)$) is a closed $n$-dimensional subset
of $T^*V=V\times V^*$ (resp. $TV$) invariant with respect to multiplication
by non-negative numbers acting on the second factor.

 Observe that ${\Nor}(K)$, and hence
$\Nor^*(A)$ as well, is biLipschitz homeomorphic to $V$: putting
$p_A:V \to A$ for the nearest point projection, $V$ maps onto
${\Nor}(A)$ via the map $ P_A: x\mapsto (p_A(x),
x-p_A(x))$, with inverse induced by $(x,y) \mapsto x+y$. It is
clear that $P_A\colon V\to V\times V$ is a proper map.

It is useful to think of these objects as defining {\it integral
currents} in the tangent and cotangent bundles of $V$. Given a
smooth manifold $M$, put  $\Omega_c^k(M)$ for the space of all
compactly supported $C^\infty$ differential forms of degree $k$ on
$M$. We recall (\cite{GMT}) that the space $\mathbb I_k(M)$ of
integral currents of dimension $k$ on $M$ is the space of all
continuous linear functionals $T: \Omega_c^k(M)\to \R$ with the
properties
\begin{itemize}
\item There is a sequence of bounded measurable subsets
$E_1,E_2,\dots \subset \R^k$ and locally Lipschitz maps $f_i :E_i \to M$ such that
$$
T = \sum_{i=1}^\infty f_{i*} \lcur E_i \rcur,
$$
where $\lcur E\rcur $ denotes the operation of integration of a
$k$-form over $E$. Note that by Rademacher's theorem the
derivative of $f$ exists almost everywhere, and constitutes a
bounded measurable function. Thus the pull-backs $f^*_i\phi$ are
integrable over the $E_i$, so the pushed forward currents
$f_{i*}\lcur E_i \rcur$, given by
$$f_{i*} \lcur E_i \rcur (\phi) := \int_{E_i} f_i^*\phi,$$
are well-defined.
\item For each compact set $C\subset M$ we have $\mass_C(T) <\infty$, where
$$ \mass_C(T):= \sup \{T(\phi): \supp \phi \subset C, \ \parallel \phi \parallel_0 \le 1\}.$$
Here $\parallel \cdot \parallel_0$ denotes the $C^0$-norm.
\item For each compact $C \subset M$,
$$
\mass_C(\partial T) <\infty.
$$
Here $\partial T$ is the current of dimension $k-1$ given by
$\partial T(\psi) := T(d\psi)$.
\end{itemize}

Equipping $M$ with a smooth Riemannian metric, the {\bf coflat
seminorm} $\parallel \phi \parallel ^\flat$ of a form $\phi \in
\Omega^k(M)$ is given by
$$
\parallel\phi \parallel^\flat:= \max\{\parallel \phi\parallel_0, \parallel d\phi\parallel_0\}.
$$
 Given $C \subset M $, we put
 $$ \parallel T\parallel^C _\flat:=
 \sup \{\left|T(\phi)\right|: \parallel \phi\parallel^\flat\le 1, \supp \phi \subset C\}.$$
This is the restriction to the lattice $\mathbb I_k(M)$ of the {\bf flat seminorm} relative to $C$.
In the case $C=M$ we will omit the superscript. The {\bf  local
flat topology} on $\mathbb I_k(M)$ is determined by the condition
that $T_1,T_2,\dots \to T$ iff
$$ \parallel T_i-T\parallel^C _\flat\to 0 $$
for every compact $C \subset M$.

{\bf Remark.} Related to the flat seminorms is the integral flat
metric
$$\mathcal F(T) := \inf\{\mass_M(R) +\mass_M(S): R\in \mathbb I_k(M), \ S \in \mathbb I_{k+1}(M), \ T = R + \partial S \}.$$
Clearly
$$
\parallel T\parallel_\flat^M \le \mathcal F(T).
$$


We now fix an orientation of $V$ and define the {\bf conic normal
cycle} of $A$ to be the integral current
\begin{equation}\label{eq:def N}
\vec N(A): = P_{A*}\lcur V\rcur \in \mathbb I_n(V \times V),
\end{equation}
where $\lcur V\rcur$ denotes the fundamental class of $V$. By
 \cite{GMT}, 4.1.14 and 4.1.24,  the image of an
integral current under a proper Lipschitz map is well-defined, and
belongs to the class of integral currents. Note that $\partial
\vec N(A) = \partial P_{A*}\lcur V\rcur = P_{A*} \partial \lcur V
\rcur = 0.$ Likewise we define the {\bf conic conormal cycle}
$\vec N^*(A)\in \mathbb I_n(V \times V^*)$ as the image of $\vec
N(A)$ under the identification $V \times V \to V \times V^*= T^*V$
arising from the euclidean structure. (Note that the image of $\vec N^*(A)$ under the antipodal map on the $V^*$ factor is identical to the {\bf characteristic cycle} of $A$.) It is easy to
see that the
supports of $\vec N(A), \vec N^*(A)$ are
${\Nor}(A),\Nor^*(A)$ respectively, and  that $\Nor^*(A)$ and $\vec N^*(A)$ are independent of the choice
of Euclidean metric.

Recall that $\K(V)$ denotes the metric space of all compact convex
subsets of $V$, endowed with the Hausdorff metric. We endow
$\mathbb I_n(V \times V^*)$ with the topology of local flat
convergence.

\begin{lemma}\label{lm:continuity}
$\vec N^*$ is continuous as a map $\K(V) \to\mathbb I_n (V
\times V^*)$.
\end{lemma}
\begin{proof} Let $A_0,A_1, \dots \in \mathcal K(V)$, with $A_i \to A_0$.
Clearly the nearest point projections $ p_{A_i}$ converge uniformly to
$p_{A_0}$, hence $P_{A_i} \to P_{A_0}$ uniformly as well.
It now follows from the discussion in 4.1.14 of \cite{GMT},
and the definition of $\vec N^*$, that $\vec N^*(A_i) \to \vec N^*(A_0)$ in the local flat topology.
\end{proof}

It is convenient to introduce a different characterization of
$\vec N^*(A)$. We say that a current in $\mathbb I_n(V\times V^*)$
is {\bf Lagrangian} if it annihilates the ideal of all multiples
of the canonical symplectic form $\omega$. The terminology is
motivated by the obvious fact that if $M \subset V \times V^*$ is
a smooth oriented $n$-dimensional submanifold then the current
given by integration of $n$-forms over $M$ is Lagrangian iff $M$
is a Lagrangian submanifold of $V \times V^*$ in the usual sense.

\begin{lemma} If $A \in \K(V)$ then $\vec N^*(A)$ is a Lagrangian current.
\end{lemma}
\begin{proof}
If $A$ has smooth boundary and nonempty interior then $\vec
N^*(A)$ decomposes as the sum of two terms. The first is
integration over the image of $A$ itself under the zero map $V \to
V^*$. This current is obviously Lagrangian. The second is
integration over the bundle of outward conormal rays to the smooth
hypersurface $\partial A$. The conormal bundle of a submanifold is
a classical example of a Lagrangian submanifold, hence the second
term is Lagrangian as well.

As is well known, every element  $A\in\K(V)$ may be approximated
in the Hausdorff metric by a sequence $A_1,A_2,\dots$ of smooth
bodies with nonempty interior. By Lemma \ref{lm:continuity}, the
conic conormal cycles $\vec N^*(A_i)$ converge in the locally flat
topology to $\vec N^*(A)$. But locally flat convergence implies
weak convergence, so the latter current must annihilate the
symplectic ideal since the the former do.
\end{proof}

We now recall the main theorem of \cite{Fu:1989}.

\begin{theorem}\label{thm:MA} Let $W$ be an oriented real vector space of dimension $n$. If $f:W \to \R$ is locally Lipschitz, then there is at most one
closed Lagrangian integral current $T\in \mathbb I_n(W \times
W^*)$ such that
\begin{itemize}
\item $T$ is locally vertically bounded, i.e. $\supp T \cap( C\times W^*)$ is compact for every compact $C \subset W$; and
\item if $\phi: W \times W^* \to \R$ is a smooth compactly supported function and $d  \vol_W$ is a (positive) volume form for $W$, then
\begin{equation}\label{eq:integral condition}
T(\phi\, \pi_W^* d\vol_W) = \int_W \phi(x,df(x)) \, d\vol_W
\end{equation}
where $\pi_W\colon W\times W^*\to W$ is the canonical projection.
\end{itemize}
If $f$ is convex then this $T$ exists.
\end{theorem}

Here $\pi_W: W \times W^*\to W$ denotes the projection to the
first factor. We will call this current $T$ the {\bf differential
cycle} of $f$, denoted here by $D(f)$.

{\bf Remarks.} 1. The differential $df(x)$ exists for almost every
$x$ by Rademacher's theorem. The resulting map $df$ is measurable,
so (\ref{eq:integral condition}) makes sense.

 2. The point is that if $f$ is smooth then the current $D(f)$ is simply integration over the graph of the differential $df$. If $f$ is convex then the graph of the subgradient
 of $f$ is a Lipschitz submanifold of $W \times W^*$ and inherits a natural orientation from that of $W$; the differential current $D(f)$ is then given by integration over the graph of the subgradient.

3. In fact a stronger form of the theorem is true: the condition
that $f$ be locally Lipschitz may be replaced by the statement $
f\in W^{1,1}_{loc}$ (i.e. $df\in L^1_{loc}$); and the first condition on $T$ may be replaced by
the requirement that the restriction $T\,  \llcorner (C\times W^*)$ have finite
mass for every compact $C \subset W$.

{\it Sketch of proof.} It is enough to show that if $T$ satisfies
the first condition, and additionally annihilates all functional
multiples of $\pi_W^*d\vol_W$, then $T = 0$. The proof of this
statement is modeled on a well-known fact about smooth Lagrangian
submanifolds $L\subset T^*W$ with the property that $\pi_W|L$ is a
submersion: locally, such a submanifold is a fiber bundle over its
image $\Lambda\subset W$, with fibers of the form $dg(x) +
\nu^*_x\Lambda$, where $g$ is a smooth function and
$\nu^*_x\Lambda$ is the conormal fiber to $\Lambda$ at $x$ (cf.
\cite{Harvey-Lawson:198?}). In particular, the fibers are
unbounded. A weak form of this description applies to a
rectifiable Lagrangian carrier of $L$ on the set of points where
the projection to $W$ has maximal rank. $\square$

Given $A\in \K(V)$, we denote by $h_A: V^* \to \R$ the support
function of $A$ given by
\begin{equation}
h_A(\lambda) := \sup_{x\in A} \lambda(x).
\end{equation}
It is well known, and easy to prove, that $h_A$ is sublinear, i.e.
convex and positively homogeneous of degree 1
\cite{schneider-book}. In particular it is Lipschitz, hence
differentiable for a.e. $y\in V^*$ by Rademacher's theorem. In
this case the differential $dh_A(\lambda)\in V$ has a particular
geometric meaning:

\begin{proposition}\label{prop:support} If $h_A$ is differentiable at $\lambda\in V^*$, then $x:=dh_A(\lambda)\in V$ is the unique point in $A$ at which $\lambda$ supports $A$, i.e. such that $\lambda(A) \subset (-\infty, \lambda(x)]$.
\end{proposition}
\begin{proof}
\cite{schneider-book},  Corollary 1.7.3.
\end{proof}

\begin{proposition}\label{prop:D&N} Let $i:V\times V^* \to V^* \times V$ denote the interchange map $i(x,y) = (y,x)$. If $A \in \K(V)$ then
\begin{equation}i_*\vec N^*(A) = D(h_A).
\end{equation}
\end{proposition}
\begin{proof} It is enough to check that $i_*\vec N^*(A)$ satisfies the conditions of Theorem \ref{thm:MA} for $f= h_A$. The first is trivial: since $\vec N^*(A)$ is supported in $A \times V^*$, it is clear that $i_*\vec N^*(A)$ is supported in $V^* \times A$ and hence is even globally vertically bounded.

To prove the second, we pass to the dual setting using our fixed
euclidean structure $( \cdot,\cdot )$ on $V$. Abusing notation,
we again denote by $h_A:V\to \R$ the support function
$$h_A(x):= \sup_{y\in A} ( x,y) .$$

Put $q_A(y) = y-p_A(y)$, and let $p_1,p_2: V\times V \to V$ be the
projections to the first and second factors, respectively.
We must show that for smooth compactly supported functions
$\varphi:V\times V\to \R$,
\begin{equation}\label{eq:axiom 2}
i_*\vec N (A) (\varphi \, p_1^*\dvol_V) = \int_V \varphi(x,\nabla
h_A(x)) \, dx.
\end{equation}
Recalling (\ref{eq:def N}) , the left-hand side may be expressed
\begin{align}
i_*P_{A*}[V](\varphi \, p_1^*\dvol_V)
&=\int_V P_A^*i^*(\varphi \, p_1^*\dvol_V) \\
&=\int_V \varphi(y-p_A(y),p_A(y)) P_A^*p_2^*\dvol_V(y)\\
&=\int_V \varphi(q_A(y),p_A(y)) q_A^*\dvol_V(y) \\
& = \int_V \varphi(q_A(y),\nabla h_A(q_A(y))) q_A^*\dvol_V(y) ,
\end{align}
by Prop. \ref{prop:support}. Since $q_A^{-1}(y)$ is a singleton
for a.e. $y\in V$, the desired relation (\ref{eq:axiom 2}) follows
from the change of variables formula.
\end{proof}

In fact the inverse map to $\vec N^*$ is also well-defined and
continuous:

\begin{corollary}\label{lm:N*inverse}
If $A,B \in \K(V)$ and $\vec N^*(A) =\vec N^*(B)$, then $A=B$. If
$A_0,A_1,A_2,\dots \in \K(V)$ and $\vec N^*(A_i) \to \vec
N^*(A_0)$ in the local flat topology, then $A_i \to A_0$ in the
Hausdorff metric topology.
\end{corollary}
\begin{proof}
If $\vec N^*(A) = \vec N^*(B)$, then $D(h_A) = D(h_B)$ by Prop.
\ref{prop:D&N}. It follows at once that $dh_A = dh_B$ a.e. in
$V^*$. Since a Lipschitz function on a euclidean space with
derivative a.e. equal to zero is constant, and $h_A(0) =h_B(0)=0$,
we conclude that $h_A = h_B$. Therefore $A=B$.

To prove the second statement, note first that all of the $A_i$
must lie within a sufficiently large fixed compact subset of $V$:
for there exists a smooth differential form $ \kappa_1$ on $T^*V$
such that $\vec N^*(A)(\kappa_1) $ is the mean breadth of $A$, for all $A
\in\K(V)$ (cf. \cite{Fu:1990}). Thus the mean breadth of the $A_i$
converges to that of $A_0$. But the mean breadth of a convex body
dominates its diameter, so we conclude that the diameters of the
$A_i$ are uniformly bounded. Furthermore there exists another
smooth differential form $\kappa_0$ such that if $\phi:V\to \R$ is
a smooth compactly supported function whose restriction to $A_0$
is equal to 1, then $\vec N^*(A_0)(\pi^*\phi \kappa_0) = 1$ where
$\pi\colon T^*V\to V$ is the canonical projection. Thus
$\vec N^*(A_i)(\pi^*\phi\kappa_0) \ne 0$ for all sufficiently large $i$,
and in particular $A_i= \pi(\supp \vec N^*(A_i))$ has a nonempty
intersection with $\supp \phi$ for such $i$.

Thus the Blaschke Selection Theorem implies that there exists a
subsequence $A_{i'}$ converging in the Hausdorff metric to some
$B_0\in \K(V)$. By Lemma \ref{lm:continuity}, $\vec N^*(A_0) =
\lim_{i'\to\infty}\vec N^*(A_{i'}) = \vec N^*(B_0)$; by the
assertion above, $A_0=B_0$. Since this outcome is independent of the chosen convergent subsequence it follows that the entire sequence of the $A_i$ converges to $A_0$.
\end{proof}

\begin{proposition}\label{prop:max} If $f,g, \min(f,g):W \to\R$ are convex, then
$$D(\max(f,g))+ D(\min(f,g)) = D(f) + D(g). $$
\end{proposition}

\begin{proof}
It is enough to show that $ D(f) + D(g)- D(\min(f,g)) $ satisfies
the conditions of Theorem \ref{thm:MA}, with $f$ replaced by the
function $\max(f,g)$. All of them are immediate except for the
last one, and for this it is enough to show that
$$ \{d(\max(f,g))(x) , d(\min(f,g))(x)\} = \{df(x),dg(x)\}$$
for a.e. $x \in W$ at which all four differentials exist (which
happens a.e. in $W$ by Rademacher's theorem).

This is obvious when $f(x) \ne g(x)$. On the other hand, let $E$
denote the set of points $x$ such that $f(x) = g(x)$ and both of
$f,g$ are differentiable at $x$. By classical measure theory, $E$
has density $1$ at a.e. point $x \in E$. If $x $ is such a point,
then clearly $df(x) = dg(x)$. Therefore this common value is also
equal to both $d\max(f,g)(x)$ and $d\min(f,g)(x)$. \end{proof}

\begin{corollary}\label{cor:N is a valuation}
 If $A,B, A\cup B \in \K(V)$ then
\begin{equation}\label{eq:sum} \vec N^*(A \cup B) + \vec N^*(A \cap B) = \vec N^*(A) + \vec N^*(B).
\end{equation}
\end{corollary}
\begin{proof}
If $A\cup B\in \K(V)$ then $h_{A\cup B} = \max(h_A,h_B)$ and
$h_{A\cap B} = \min(h_A,h_B)$. Therefore (\ref{eq:sum}) follows at
once from Props. \ref{prop:D&N} and \ref{prop:max}.
\end{proof}

It is sometimes convenient to consider instead the (non-conic)
{\bf normal cycle} $N(A)$ in the tangent sphere bundle $V\times
S(V)$, and the corresponding {\bf conormal cycle} $N^*(A)$ in the
cotangent ray bundle $V \times \PP_+(V^*)$. To define $N(A)$, let
$r:V\to [0,\infty)$ denote the length function induced by the
fixed euclidean metric. Then
\begin{equation}\label{eq:N}
N(A):= g_*\left(\langle \vec N(A), r\circ p_2, 1\rangle\right) \in
\mathbb I_{n-1}(V \times S(V)),
\end{equation}
where $\langle T,f ,c\rangle $ denotes the slice of the current
$T$ by the function $f$ at the value $c$ (cf. \cite{GMT}, 4.3) and
$g: V \times (V-\{0\}) \to V\times S(V)$ is the normalizing map
$(x,y) \mapsto (x, \frac{y}{|y|})$. Note that the slicing
operation is well-defined for a.e. $c$ whenever $T$ is an integral
current and $f$ is a Lipschitz function whose restriction to the
support of $T$ is proper, and may be thought of as the
intersection of $T$ with the level set $f^{-1}(c)$.  In the
present case the slice is well-defined at {\it every} value of $r$
since, putting $\theta_c(x,y) := (x,cy)$ for $c >0$,
\begin{equation*}
\langle \vec N(A), r, a\rangle = \theta_{c*}\langle \vec N(A), r,
\frac a c\rangle
\end{equation*}
--- this in view of the facts
\begin{align*}
\theta_{c*} \vec N(A) &= \vec N(A),\\
r\circ \theta_c  &= c\, r
\end{align*}
and the general formula
$$ \langle h_*T, f, c\rangle = h_* \langle T, f\circ h, c\rangle $$
(cf. \cite{GMT}, 4.3.2(7)).

The conormal cycle $N^*(A)$ is then the image of $N(A)$ under the
natural map $V \times S(V)\to V\times \PP_+(V^*)$ induced by the
euclidean metric. It is clear that $N^*(A)$ does not depend on the
choice of this metric. Recall that $V\times\PP_+(V^*)$ carries a
natural contact structure, and a choice of metric even determines
a particular global contact 1-form. The current $N^*(A)$ is
Legendrian in the sense that it annihilates every element of the
ideal generated by any such 1-form.

The conic normal cycle may be reconstructed from the normal cycle
in a canonical way.  Given a manifold $M$, put $\mathbb I_k^c(M)$
for the space of all compactly supported integral currents of
dimension $k$ on $M$, where the topology on this space is
determined by the condition that $T_i\to T_0\in \mathbb I^c_k(M)$
iff $T_i\to T_0$ in the local flat topology and all of the $T_i$
are supported in a single compact set $C\subset M$. Define $f: V
\times S(V)\times \RR \to V \times V$ by $f(x,v;t):= (x, tv)$.
Define also $g: V \times \R\to V$ by $ g (x, t): = tx$ and put $h:
V \times S(V) \to V$ for the projection.  Now define $\gamma:
\mathbb I_{n-1}^c(V\times S(V)) \to\mathbb I_n( V \times V)$ by
\begin{equation}\label{eq:N to vec N}
\gamma(T):= f_*\left(T \times [0,\infty)  \right) + g_*(h_*T
\times \lcur 0,1\rcur) \times \lcur 0\rcur .
\end{equation}
Since $f$ is proper this map is continuous in view of the topology
given above on $\mathbb I^c_{n-1}$. If $\partial T = 0$ then the
first factor in the second term above may be characterized as the
unique compactly supported current in $V$ with boundary equal to
$h_*T$.

\begin{proposition}\label{prop:N to vec N} If $A\in \K(V)$ then $\gamma(N(A)) = \vec N(A)$.
\end{proposition}
\begin{proof} For $r\ge 0$, put $A_r:= \{x\in V: \dist (x,A) = r\}$. Then $N(A) = P_{A*}\lcur A_1\rcur$ and
 the map $Q:(x,t) \mapsto t(x-p_A(x)) + p_A(x) $ is an
 orientation-preserving locally biLipschitz homeomorphism $A_1\times (0,\infty) \to V\setminus A$.
 These maps satisfy the relations
\begin{align*}
P_A \circ Q&= f\circ (P_A \times \operatorname{id}), \\
h\circ P_A &= p_A.
\end{align*}
Therefore
$$ h_*N(A) = p_{A*}\lcur A_1 \rcur  =\partial \lcur A\rcur,$$
and, using the characterization above of the second term in
(\ref{eq:N to vec N}),
\begin{align*}
\vec N(A) &:=P_{A*} \lcur V\rcur\\
&=P_{A*} \lcur A \rcur + P_{A*}\lcur V\setminus A\rcur \\
&=\lcur A\rcur \times \lcur 0\rcur + P_{A*} Q_* \lcur A_1\times (0,\infty) \rcur \\
&=g_*(\partial \lcur A\rcur \times \lcur 0,1\rcur) \times \lcur 0\rcur+
f_*\left(P_{A*} \lcur A_1 \rcur \times \lcur 0,\infty \rcur\right)\\
&= g_*(h_*N(A)\times \lcur 0,1\rcur) \times \lcur 0\rcur +
f_*\left( N(A) \times \lcur 0,\infty \rcur\right),
\end{align*}
as claimed.
\end{proof}

It is immediate from the definition that the normal and conormal
cycles share the basic properties of their conic counterparts
described above. For brevity we give the explicit statements
only in the conormal case:

\begin{proposition}\label{prop:conormal}
$N^*$ is a continuous injection $\K(V) \to \mathbb I_{n-1}( V
\times \PP_+( V^*))$, and the inverse map (defined on the image)
is continuous, where the topology on $\mathbb I_{n-1} ( V \times
\PP_+( V^*))$ is the local flat topology. $N^*$ is a valuation in
the sense that if $A,B,A\cup B \in \K(V)$, then $N^*(A\cup B) +
N^*(A\cap B) = N^*(A) + N^*(B)$.
\end{proposition}

Put $C^1_\flat $ for the normed space of all $C^1$-smooth
differential forms $\phi$ of degree $n-1$ on $ (V \times \PP_+(
V^*))$, with finite coflat norm.
 The preceding discussion yields the following.
\begin{theorem} The map $\K(V) \times C^1_\flat \to \R$ given by
$$ (A, \phi) \mapsto N^*(A)(\phi)$$
is continuous.
\end{theorem}
{\bf Remark.} The space $C^1_\flat $ may be replaced by the space
of all flat cochains (cf. \cite{Whitney:1957}, p. 233) with finite
coflat norm .
\begin{proof} Suppose $A_1,A_2,\dots \to A_0$ in $\K(V)$ and
$\phi_1,\phi_2,\dots\to \phi_0$ in $C^1_\flat$.
By Prop. \ref{prop:conormal}, given $\eps >0$ there is $M \in \mathbb N$ such that
\begin{align*}
\parallel N^*(A_i)-N^*(A_0))\parallel_\flat &< \eps, \\
\parallel \phi_i-\phi_0 \parallel^\flat &< \eps \\
\end{align*}
for $i \ge N$. Thus
\begin{align*}|N^*(A_i)(\phi_i) - N^*(A_0)
(\phi_0 ) |& \le |N^*(A_i)(\phi_i) - N^*(A_0) (\phi_i )| + |N^*(A_0)(\phi_i) - N^*(A_0) (\phi_0 ) |  \\
&\le \parallel N^*(A_i)-N^*(A_0) \parallel_\flat
\parallel \phi_i\parallel^\flat + \mass(N^*(A_0))\parallel \phi_i-\phi_0\parallel_0 \\
&< \eps( \parallel \phi_0 \parallel^\flat + \, \eps + \mass
N^*(A_0)),
\end{align*}
which proves the desired assertion.
\end{proof}

\subsection{Normal cycles for more general sets.}
Normal cycles are also available for a wide class of sets other
than convex ones. First let us define it in the class of sets
presentable as finite unions of convex compact sets. Let
\begin{eqnarray}\label{i1}
X=\cup_{i=1}^N A_i,\, A_i\in \ck(V).
\end{eqnarray}
Set
\begin{eqnarray*}
N^*(X):=\sum_{I\subset\{1,\dots,N\}, I\ne \emptyset}
(-1)^{|I|+1}N^*(\cap_{i\in I}A_i),\\
\vec N^*(X):=\sum_{I\subset\{1,\dots,N\}, I\ne \emptyset}
(-1)^{|I|+1}\vec N^*(\cap_{i\in I}A_i).
\end{eqnarray*}
Using Corollary \ref{cor:N is a valuation} it is easy to check
that the definitions of $N^*(X)$ and $\vec N^*(X)$ do not depend on
a choice of presentation (\ref{i1}).

Let us define $N^*(X)$ and $\vec N^*(X)$ where $X$ is a compact
smooth submanifold with boundary. For any point $x\in X$ let us
define the {\itshape tangent cone} to $X$ at $x$, denoted by $T_xX$,
the set
$$T_xX:=\{\xi\in T_xV| \mbox{ there exists a } C^1-\mbox{map }
\gamma\colon [0,1]\to X \mbox{ such that }\gamma'(0)=\xi\}.$$ It
is easy to see that $T_xX$ coincides with the usual tangent space
if $x$ is an interior point of $X$, and it is a halfspace if $x$
belongs to the boundary of $X$. Define
\begin{eqnarray}\label{i2}
\Nor^*(X):=\cup_{x\in X}-(T_xX)^o
\end{eqnarray}
where for a convex cone $C$ in a linear space $W$ one denotes
$C^o$ its dual cone in $W^*$:
$$C^o:=\{y\in W^*|\, y(x)\geq 0\, \mbox{ for any } x\in C\}.$$
Clearly $\Nor^*(X)$ is invariant under the group $\RR_{>0}$ of
positive real numbers acting on the cotangent bundle $T^*V$ by
multiplication along the fibers.
The sets $\Nor(X)$ and $\Nor_1(X)$ are now defined as in Definition \ref{cc-1-2}.

The corresponding current $N^*(X)$ is more naturally constructed in the more general context of semi-convex sets  $X$ (cf. \cite{zahle90}, \cite{Federer:1959}). Put $p_X:X_{[0,r)} \to X$ for the nearest point projection to $X$, defined for the set $X_{[0,r)}:= \{x \in V: \dist(x,X) < r\}$, where $r:=\reach(X)$ of $X$, and $P_X(x):= (p_X(x), x-p_X(x))$. Put
$$ \vec N'(X): = P_{X*}\lcur X_{[0,r)} \rcur,$$
where the domain $X_{[0,r)}$ inherits its orientation from $V$. Choose a diffeomorphism $f:[0,r)\to [0,\infty)$, and put $F(x,v):= (x, f(|v|) \frac{v}{|v|})$. Now put
$$ \vec N(X):= F_*\vec N'(X).$$
>From this current we may construct the currents $\vec N^*(X), N(X), N^*(X)$ as in the remarks surrounding (\ref{eq:N}).
 It is easy to see that if $X$ is convex then
these definitions of $N^*(X)$ and $\vec N^*(X)$ coincide with the
previous ones. Furthermore $\supp N^*(X) = \Nor^*(X)$ and $\supp N(X) = \Nor_1(X)$.

The normal cycle of subanalytic subsets was defined in \cite{fu:1994} using tools from
geometric measure theory (in fact Thm. 3.2 of \cite{fu:1994} gives a unique characterization of the normal cycle of an arbitrary compact set in $\mathbb R^n$, dual to Thm. \ref{thm:MA}). A similar notion of chacteristic cycle of subanalytic subsets was
introduced independently in \cite{kashiwara-schapira}
using tools from sheaf theory.

Conormal cycles transform in a natural way under diffeomorphisms. We will only need this fact in the smooth case:

\begin{lemma}\label{lift lemma} Let $X\subset V$ be a compact domain with smooth boundary, and $U \supset X$ an open neighborhood. Let $f:U \to W\subset V$ be an orientation-preserving diffeomorphism, and put $\bar f: U \times \mathbb P_+(V^*) \to W \times \mathbb P_+( V^*)$ be the natural lift of $f$ defined by
$$\bar f (x, [\lambda]) :=  (f(x), [(f^{-1})^*\lambda]).$$
Then
$$ \bar f_* N^*(X) = N^*(f(X)).$$
\end{lemma}
\begin{proof}
It is clear that $\bar f$ maps the manifold of outward conormal lines to $\partial X$ diffeomorphically onto that of $f(\partial X) = \partial f(X)$. The cycles $N^*(X)$ and $N^*(f(X))$ are the fundamental classes of these manifolds, and therefore $\bar f_*N^*(X) = \pm N^*(f(X))$. To see that the sign is positive, we note that the orientations of the fundamental classes are determined by  the relations
$$ p_{1*} N^*(X)= \partial \lcur X\rcur, \  p_{1*} N^*(f(X))= \partial \lcur f(X)\rcur.$$
Since $f$ preserves orientation by hypothesis,
$$p_{1*} \bar f_*N^*(X) = f_*\partial \lcur X \rcur = \partial \lcur f(X) \rcur = p_{1*}N^*(f(X)),$$
which establishes the claim.
\end{proof}

M. Z\"ahle has proved the following fundamental approximation
theorem:

\begin{theorem} \label{thm:poly approx}Let $X \subset V$ be a compact domain with smooth boundary.
There exists a sequence of polyhedra $P_1,P_2,\dots\subset V$ such that
$$\lim_{i\to \infty} N(P_i) = N(X).$$
\end{theorem}
\begin{proof} This is \cite{zahle90}, Theorem 1. The proof given there may be simplified as follows.

One may show, along the lines of \cite{zahle90} or \cite{Fu:1993},
that there exists a sequence of polyhedra $P_i$ and a constant $C<
\infty$ such that
\begin{align}
\label{eq:one} \mass(N(P_i)) &\le C, \\
\label{eq:two} \supp N(P_i) &\to \supp N(X),\\
\label{eq:three} \lcur P_i\rcur& \to \lcur X\rcur.
\end{align}
Clearly the relation (\ref{eq:three}) implies that $\partial
\lcur P_i\rcur \to \partial \lcur X\rcur$.

By the compactness theorem for integral currents (\cite{White})
and the constancy theorem (\cite{GMT}, p. 357), the relations
(\ref{eq:one}) and (\ref{eq:two}) imply that there is a
subsequence $N(P_{i'}) \to k N(X)$ for some integer $k$. Thus
$\partial\lcur P_{i'}\rcur = \pi_*N(P_{i'}) \to k\pi_* N(X)=
k\partial \lcur X\rcur$, so (\ref{eq:three}) implies that $k =1$.
Since this is independent of the subsequence chosen, the result
follows.
\end{proof}

\begin{corollary} \label{cor:vec approx} In the scenario of Theorem \ref{thm:poly approx} we have also
$$ \lim_{i\to \infty} \vec N(P_i) = \vec N(X)$$
\end{corollary}
\begin{proof} This follows at once from Theorem \ref{thm:poly approx} and Prop. \ref{prop:N to vec N}
since the map $\gamma$ occurring there is continuous.
\end{proof}

\section{Auxiliary results.}\label{auxiliary} \setcounter{subsection}{1}
\setcounter{theorem}{0} The goal of this section is to prove some
technical results which will be used in the construction of the
product on valuations in Section \ref{construction}. The main
results of this section are Lemmas \ref{l8.05}, \ref{l9}, and
Proposition \ref {l10}.

Let $V$ be an $n$-dimensional real vector space. As usual we fix a Euclidean metric on $V$. In this section we will also fix
a compact smooth $n$-dimensional submanifold with boundary
$X'\subset V\times V$ which projects diffeomorphically onto its
images in $V$ under both projections $p_1,p_2\colon V\times V\to
V$. We also fix a compact submanifold with boundary $X\subset X'$
such that $X\cap \pt X'=\emptyset$. We will denote throughout this
section $$\tilde p_1,\tilde p_2\colon X'\to V$$ the restrictions
of the projections $p_1,p_2$ respectively to $X'$.

For a domain $\Ome\subset V$ with smooth boundary let us denote by
$k_1(\Ome,s),\dots,k_{n-1}(\Ome,s)$ the principal curvatures at a
point $s\in \pt\Ome$.
\begin{lemma}\label{l2}
Let $A_1,\dots,A_k\subset V$ be compact strictly convex subsets
with smooth boundaries. Then there exists a constant $C>0$
(depending on these subsets) such that for any
$\lam_1,\dots,\lam_k\geq 0$ with $\sum_{i=1}^k\lam_i=1$ one has
for any $1\leq l\leq n-1$
$$\frac{1}{C}\leq k_l(\sum_{i=1}^k\lam_iA_i,s)\leq C$$
for any $s\in \pt(\sum_{i=1}^k\lam_iA_i)$.
\end{lemma}
{\bf Proof.} For convenience let us fix on $V$ an orthonormal
coordinate system $(x_1,\dots,x_n)$.The principal curvatures $k_l(A)$ of a strictly convex compact set $A$ can be estimated from both
sides via the eigenvalues of the Hessian of the supporting
functional $Hess(h_A):=\left(\frac{\pt^2h_A}{\pt x_i\pt
x_j}\right)$ restricted to the tangent bundle to the unit sphere
$TS^{n-1}$. Let us denote by $H_i:=Hess(h_{A_i})|_{TS^{n-1}}$,
 $i=1,\dots,k$.

Let us denote by $c'$ and $C'$ the minimum and the maximum
respectively over the unit sphere $S^{n-1}$ of all the eigenvalues
of all $H_i$'s. Then clearly $0<c'<C'<\infty$. Then
$$(\sum_i\lam_i)c'Id_{n-1}\leq \sum_{i}\lam_iH_i\leq
(\sum_i\lam_i)C'\cdot Id_{n-1}.$$ The lemma follows. \qed


\begin{proposition}\label{l3}
Let $\delta_0,\delta_1 >0$ be given. Suppose $X\subset V$ is a
compact subset with $\reach(X)>\delta_0$. Suppose also that
$A\subset V$ is a strictly convex compact set with smooth
boundary, containing the origin in the interior, and such that all
principal curvatures $k_l$ of $\pt A$ satisfy $k_l \ge \delta_1$.
Then for any $0<\eps \le \delta_0\delta_1$, the map
$$\xi_{A,\eps} \colon V\times \PP_+(V^*)\times [0,1]\to V$$
given by $$\xi_{A,\eps}(p,n,t)= p+\eps \cdot t\nabla h_{A}(n)$$ is
a homeomorphism of $N(X)\times (0,1]$ onto  $(X+\eps A)\backslash
X$.
\end{proposition}

For the proof we will need the two assertions, both known, of the following lemma.
\begin{lemma}\begin{itemize}
\item If $\reach (X) > \delta_0$ and $n_i \in \Nor(X,x_i), i = 0,1$, with $ |n_0| = |n_1| = 1$, then
\begin{equation}\label{federer's inequality}
\left(x_1-x_0 ,n_1-n_0\right) \ge - \delta_0^{-1} |x_1-x_0|^2.
\end{equation}
\item Let $A$ be a convex body with smooth boundary, with all principal curvatures $k_l \ge \delta_1$. Suppose
$x_i \in \pt A$, with outward normals $n_i$, $ i = 0,1$. Then
\begin{equation}\label{convex inequality}
\left(x_1-x_0, n_1-n_0\right)\ge \delta_1 |x_1-x_0|^2.
\end{equation}
\end{itemize}
\end{lemma}
\begin{proof} The first assertion follows at once from \cite{Federer:1959}, Theorem 4.8 (7).

The second assertion may be deduced as follows. It is easy to see using Schur's theorem (\cite{Chern}) that if
$B$ is a ball of radius $\delta_1^{-1}$ passing through $x_0$, and with outward normal $n_0$ there, then $B \supset A$.
On the other hand an elementary calculation shows that
$(n_0 , x_0-p) \ge \frac{\delta_1} 2 |p-x_0|^2$ for every $p \in B$.
In particular  $(n_0 , x_0-x_1) \ge \frac{\delta_1} 2 |x_1-x_0|^2$.
Adding this to the analogous inequality for $x_1$ gives (\ref{convex inequality}).
\end{proof}

\begin{proof} (of Proposition \ref{l3}) The proposition is equivalent to the assertion that if $r \le\delta_0\delta_1$
then there is no translate $p- rA$ of $r(-A)$ with interior disjoint from $X$ and intersecting $X$
in two distinct points.

Suppose there is such a translate, with $x_0,x_1\in X \cap
(p-rA)$, $ x_0 \ne x_1$. Let $a_i:= r^{-1}(p-x_i)\in \pt A$, and let $n_i$ be the outward unit normals to
$A$ at $a_i$, $i=0,1$. By (\ref{convex inequality}),
\begin{equation}
r^{-1}(x_0-x_1,n_1-n_0)=(a_1-a_0, n_1-n_0) \ge \delta_1|a_1-a_0|^2 = \delta_1
r^{-2}|x_1-x_0|^2,
\end{equation}
or
\begin{equation}
(x_0-x_1,n_1-n_0) \ge r^{-1}\delta_1|x_1-x_0|^2.
\end{equation}
On the other hand it is clear that the $n_i \in \Nor(X,x_i)$, so (\ref{federer's inequality}) gives
\begin{equation}
(x_0-x_1,n_1-n_0) \le \delta_0^{-1} |x_0-x_1|^2.
\end{equation}
Thus  $r \ge \delta_0\delta_1$, as claimed.
\end{proof}


Recall that we denote by $L$ the (real) line bundle over
$\PP_+(V^*)$ whose fiber over $l\in \PP_+(V^*)$ is equal to the
space of $\RR$-valued linear functionals on $l$. Let us denote by
$\ch:=C^\infty(\PP_+(V^*),L)$ the Fr\'echet space of
$C^\infty$-smooth sections of $L$. Clearly $\ch$ coincides with
the space of smooth functions on $V^*\backslash \{0\}$ which are
homogeneous of degree one.
\begin{lemma}\label{l4}
Consider the maps
\begin{eqnarray*}
\xi\colon \ch\times V\times \PP_+(V^*)\times[0,1]\to V
\end{eqnarray*}
defined by
\begin{eqnarray}\label{b0.1}
\xi(h,p,n,t)=p+t\cdot \nabla h(n),
\end{eqnarray} and
\begin{eqnarray*}
\Xi\colon \ch\times C^\infty (V,|\ome_V|)\to C^\infty(V\times
\PP_+(V^*)\times [0,1],\Ome^n\otimes p^*o)
\end{eqnarray*}
defined by
$$\Xi(h,\eta)=(\xi(h,\cdot))^* \eta$$
where $\xi(h,\cdot)\colon V\times \PP_+(V^*)\times[0,1]\to V$ is
defined by (\ref{b0.1}).

Then $\Xi$ is an infinitely smooth map (of infinite dimensional
manifolds), and it is linear with respect to the second argument.
\end{lemma}
{\bf Proof} is obvious. \qed

For $k\in \ZZ_+$ let us denote by $\Xi_k$, $\frac{1}{k!}$ times
the $k$-th differential of $\Xi$ at 0 with respect to the first
argument. Namely for $h_1,\dots,h_k\in \ch,\eta \in
C^\infty(V,|\ome|)$
$$\Xi_k(h_1,\dots,h_k,\eta)=\frac{\pt^k}{\pt\lam_1\dots
\pt\lam_k}\big |_0
\left(\xi(\sum_{i=1}^k\lam_ih_i)\right)^*\eta.$$ Thus
$$\Xi_k\colon \ch^k\times C^\infty(V,|\ome_V|)\to C^\infty(V\times
\PP_+(V^*)\times[0,1],\Ome^n\otimes p^*o)$$ is a continuous map
linear with respect to all $k+1$ arguments.

By the L. Schwartz kernel theorem, $\Xi_k$ extends canonically to a
continuous linear operator
\begin{eqnarray}\label{b1}
\Xi_k\colon C^\infty(V\times (\PP_+(V^*))^k, |\ome_V|\boxtimes
L^{\boxtimes k})\to C^\infty (V\times
\PP_+(V^*)\times[0,1],\Ome^n\otimes p^*o).
\end{eqnarray}
(Note that we denote this operator by the same letter $\Xi_k$.)

Let us also denote
\begin{eqnarray}\label{b2}
\Theta\colon C^\infty(V\times \PP_+(V^*)\times[0,1],\Ome^n\otimes
p^*o)\to V^\infty(V)
\end{eqnarray}
the canonical map given by integration with respect to the
normal cycle times the segment $[0,1]$. Namely
$(\Theta(\ome))(P)=\int_{N^*(P)\times [0,1]}\ome$ for any $P\in
\cp(V)$.
\begin{lemma}\label{l4.1}
Let $\psi\in V^\infty(V)$.  Assume that there exist $k\in \NN$,
sequences of smooth densities $\{\mu_N\}_{N=1}^\infty\subset
C^\infty(V,|\ome_V|)$, and $\{B_N^i\}_{N=1}^\infty\subset
\ck(V),\, i=1,\dots,k,$ of strictly convex compact sets with
smooth boundaries, containing the origin in the interior, and such
that

1) for any compact subset $T\subset V$ and any $L\in \NN$
\begin{eqnarray}\label{b3}
\sum_{N=1}^\infty ||\mu_N||_{C^L(T)}\cdot \prod_{i=1}^k
||h_{B_N^{i}}||_{C^L(S^{n-1})}< \infty;\end{eqnarray}

2) for any $K\in \ck(V)\cap \cp(V)$ one has
\begin{eqnarray}\label{b4}
\psi(K)=\sum_{N=1}^\infty \frac{\pt^k}{\pt\lam_1\dots \pt
\lam_k}\big |_0 \mu_N(K+\sum_{i=1}^k\lam_i B_N^i).
\end{eqnarray}

Then one has
\begin{eqnarray}\label{b5}
\psi=\sum_{N=1}^\infty (\Theta\circ \Xi_k)\left(\mu_N\otimes
h_{B_N^1}\otimes \dots\otimes h_{B_N^k}\right)
\end{eqnarray}
where the last series converges in the space $V^\infty(V)$.
\end{lemma}
{\bf Proof.} The inequality (\ref{b3}) implies that the series
$\sum_{N=1}^\infty \mu_N\otimes h_{B_N^1}\otimes \dots\otimes
h_{B_N^k}$ converges in $C^\infty(V\times (\PP_+(V^*))^k,
|\ome_V|\boxtimes L^{\boxtimes k})$. Hence the series in the right
hand side of (\ref{b5}) converges in $V^\infty(V)$ due to the
continuity of $\Theta$ and $\Xi_k$. Let us denote its limit by
$\psi'$.

By Lemma 2.4.5 of \cite{part2} any smooth valuation is defined
uniquely by its values on sets from $\ck(V)\cap \cp(V)$. Hence it
is enough to check that for any $K\in \ck(V)\cap \cp(V)$ one has
$$\psi(K)=\psi'(K).$$
By continuity we may assume that
$$\psi(K)=\frac{\pt^k}{\pt\lam_1\dots \pt \lam_k}\big |_0
\mu(K+\sum_{i=1}^k\lam_i B^i)$$ for any $K\in \ck(V)\cap \cp(V)$,
and thus $\psi'=(\Theta \circ \Xi_k)(\mu\otimes h_{B^1}\otimes
\dots \otimes h_{B^k})$.

Fix $K\in \ck(V)\cap \cp(V)$. By Lemma \ref{l2} and Proposition
\ref{l3} the map $\Nor_1(K)\times (0,1]\to V$ given by
$$(p,n,t)\mapsto p+t\sum_{i=1}^k\lam_i\nabla h_{B^i}(n)$$
is a homeomorphism of $\Nor_1(K)\times (0,1]$ onto
$(K+\sum_{i=1}^k\lam_i B^i)\backslash K$ for $0<
\lam_1,\dots,\lam_k\ll 1$. Hence
\begin{eqnarray*}
\mu(K+\sum_{i=1}^k\lam_i B^i)=\mu(K) +\int_{N^*(K)\times (0,1]}
(\xi(\sum_{i=1}^k\lam_i \nabla h_{B^i}))^*\mu= \mu(K)+
\int_{N^*(K)\times (0,1]} \Xi(\sum _{i=1}^k\lam_ih_{B^i},\mu).
\end{eqnarray*}
Hence
\begin{eqnarray*}
\psi(K)=\int_{N^*(K)\times[0,1]} \Xi_k(h_{B^1},\dots,h_{B^k},\mu)=
\left((\Theta\circ\Xi_k)(\mu\otimes h_{B^1}\otimes \dots \otimes
h_{B^k})\right)(K)=\psi'(K).
\end{eqnarray*}
\qed


We will also need the following simple lemma.
\begin{lemma}\label{l5}
For any subsets $T\subset V\times V$, $A,B\subset V$, and any
$x_0\in V$ one has
$$(T+(A\times B))\cap (\{x_0\}\times V)=\{x_0\}\times
(p_2(T\cap p_1^{-1}(x_0-A))+B).$$
\end{lemma}
{\bf Proof.} We have
$$(T+(A\times B))\cap (\{x_0\}\times V)=
\left((T+(A\times\{0\}))\cap (\{x_0\}\times
V)\right)+(\{x_0\}\times B).$$ Then we have for any $y\in V$
\begin{eqnarray*}
(x_0,y)\in (T+(A\times\{0\}))\cap (\{x_0\}\times V)
\Leftrightarrow \\
 \exists x\in V  \exists a\in A \mbox{ s. t. }
x+a=x_0 \mbox{ and } (x,y)\in T \Leftrightarrow \\
\exists x\in V \mbox{ s.t. } (x,y)\in T\cap p_1^{-1}(x_0 -A).
\end{eqnarray*}
The result follows. \qed

In general we will denote by $||\cdot ||_0$ the $C_0$-norm of a
map.

\begin{lemma}\label{l6}
Let $\cd,\cd'\subset V$ be compact domains with smooth boundaries,
and $g\colon \cd\tilde\to \cd'$ be a diffeomorphism. Let $K
\subset \cd$ be a compact convex subset, and put $\delta$ for the
distance from $K$ to the complement of $\cd$. Then
\begin{equation}\label{eq:reach ineq}
\reach g(K)\ge \min \left\{ \frac \delta 2 \parallel
D(g^{-1})\parallel_0^{-1}, \parallel
D(g^{-1})\parallel_0^{-2}\,\parallel D^2g\parallel_0^{-1}\right\}.
\end{equation}
\end{lemma}
{\bf Proof.} Suppose $\reach g(K) <\frac{\delta}{2} \parallel D(g^{-1})\parallel_0^{-1}$. Let $\eps >0$ be given.
 Then there exist points $q \in V$ and $x_0,x_1 \in K$, $ x_0 \ne x_1$, such that
\begin{align*}  |q- g(x_i)| &= \dist(q,g(K)) < \frac \delta 2 \parallel
 D(g^{-1})\parallel_0^{-1}, \\
 |q- g(x_i)| & < \reach g(K) + \eps, \ \  i = 0,1.
 \end{align*}
 We may assume for simplicity that $x_0 = g(x_0) = 0$.
The mean value theorem implies that the distance from $g(K)$ to
the complement of $\cd'$ is at least $\delta \, \parallel
D(g^{-1}) \parallel_0^{-1} > 2 |q|$, so by the triangle inequality
the line segment joining $ 0 $ to $g(x_1)$ lies in $\cd'$.

Since $0$ and $g(x_1)$ both lie on the sphere of radius $|q|$
about $q$, it follows that
\begin{equation}\label{eq:ineq}
(q , g(x_1)) = \frac{|g(x_1)|^2}{2}  \ge \frac{|x_1|^2}{2}  \,
\parallel D(g^{-1}) \parallel_0^{-2}
\end{equation}
by the mean value theorem. Abbreviating $L:= Dg(0)$ and letting
$\Hat L$ be its adjoint, it is clear that $\Hat L(q) \in
\Nor(K,0)$. Thus by (\ref{eq:ineq}),
\begin{align*}
\frac{|x_1|^2}{2}  \,\parallel D(g^{-1}) \parallel_0^{-2} &\le
    (L(x_1) , q) + (g(x_1) - L(x_1), q )\\
& = (x_1 , \Hat L    (q)) + (g(x_1) - L(x_1), q) \\
&\le (g(x_1) - L(x_1), q)\\
&\le \frac{|x_1|^2}{2}\,\parallel D^2g\parallel_0 |q|
\end{align*}
by Taylor's theorem and the Cauchy-Schwartz inequality. Thus
$$\reach g(K) +\eps > |q| \ge \parallel D(g^{-1}) \parallel_0^{-2} \,  \parallel   D^2 g\parallel_0^{-1} .$$
\qed


\begin{lemma}\label{l7}
Let $X\subset X'\subset V\times V$ be submanifolds as in the
beginning of this section. Assume moreover that $p_1(X)$ is
convex.

Then there exists $\delta >0$ (depending on the $C^2$-norm of the
map $(p_2\circ p_1^{-1})|_{p_1(X')}$ and its inverse, and the
distance from $X$ to $\pt X'$ only) such that for any $A\in
\ck(V)$ and any $x\in V$ the set $(X+(A\times \{0\}))\cap
(\{x\}\times V)$ is either empty or has reach at least $\delta$ as
a subset of $\{x\}\times V$.
\end{lemma}
{\bf Proof.} Lemma \ref{l5} and the assumptions imply that
\begin{equation}\label{intersection}
(X+(A\times \{0\}))\cap \{x\}\times V=\{x\}\times (\tilde p_2\circ
\tilde p_1^{-1})(p_1(X)\cap (x-A)).
\end{equation}
Now the proof follows
immediately from Lemma \ref{l6}. \qed



>From Lemmas \ref{l5} and \ref{l7} we deduce immediately the
following corollary.
\begin{corollary}\label{l8}
Let $X\subset X'\subset V\times V$ be as at the beginning of this
section. Assume that $p_1(X)$ is convex.

Then for any $x\in V$ and any $\psi\in V^\infty(V)$ one has
$$\psi((X+(A\times \{0\}))\cap (\{x\}\times V))=\left((\tilde p_1\circ \tilde
p_2^{-1})_*\psi\right)(p_1(X)\cap (x-A))$$ where $\tilde
p_1,\tilde p_2\colon X'\to V$ are the restrictions of the
projections $p_1,p_2$ to $X'$.
\end{corollary}
\begin{lemma}\label{l8.05}

(1) The function
$$\gamma\colon V^\infty(V)\times C^\infty(V,|\ome_V|)\times
 \ck(V)^{k+1}\times \RR_{\geq 0}^k\to \CC$$ defined by
$$(\phi;\mu;K,A^1,\dots,A^k;\lam_1,\dots,\lam_k)\overset{\gamma}{\mapsto} \int_{x\in
V}\phi(K\cap(x-\sum_{i=1}^k\lam_iA^i))d\mu(x)$$ with $\phi\in
V^\infty(V),\mu\in C^\infty(V,|\ome|), (K,A^1,\dots,A^k)\in
\ck(V)^{k+1},\, \lam_i\geq 0$, is a continuous function which is
$C^\infty$-smooth on $V^\infty(V)\times C^\infty(V,|\ome|)\times
\RR_{\geq 0}^k$ for fixed $(K,A^1,\dots,A^k)\in \ck(V)^{k+1}$.

(2) Fix $R>0,\, k\in \NN$. Then there exist a constant $C$, a
positive integer $L\in \NN$, and continuous seminorms $||\cdot ||$
and $||\cdot||'$ on $V^\infty(V)$ and $C^\infty(V,|\ome_V|)$
respectively depending on $n,\, k,$ and $R$ only, such that for
any strictly convex compact sets $A^1,\dots,A^k$ with smooth
boundaries, and any $K\in \ck(V)$ such that $K$ is contained in
the centered Euclidean ball of radius $R$, one has an estimate
\begin{eqnarray}\label{c1}
\big |\frac{\pt^k}{\pt\lam_1\dots \pt\lam_k}\big |_0\int_{x\in V}
\phi(K\cap(x-\sum_{i=1}^k\lam_iA^i))d\mu(x)\big |\leq
||\phi||\cdot
||\mu||'\cdot\prod_{i=1}^k||h_{A^i}||_{C^L(S^{n-1})}.
\end{eqnarray}
\end{lemma}
{\bf Proof.} By Theorem \ref{onto} there exists a continuous
epimorphism of Fr\'echet spaces
\begin{eqnarray*}
\ct\colon \oplus_{k=0}^nC^\infty(V\times
\PP_+(V^*)^k,|\ome_V|\boxtimes L^{\boxtimes k}) \surj SV(V)
(\tilde \leftarrow V^\infty(V))
\end{eqnarray*}
which is uniquely characterized by the property that for any $K\in
\ck(V)$, any strictly convex compact sets with smooth boundary
$A_1,\dots,A_k$, and any $\nu\in C^\infty(V,|\ome_V|)$ one has
\begin{eqnarray*}
\ct(\nu\boxtimes h_{A_1}\boxtimes \dots\boxtimes h_{A_k})(K)=
\frac{\pt^k}{\pt\lam_1\dots \pt\lam_k}\big|_0\nu(K+\sum_{i=1}^k\lam_iA_i).
\end{eqnarray*}
By the Banach inverse mapping theorem the map $\ct$ induces an isomorphism of Fr\'echet
spaces
\begin{eqnarray*}
\left(\oplus_{k=0}^nC^\infty(V\times
\PP_+(V^*)^k,|\ome_V|\boxtimes L^{\boxtimes k})\right)/Ker \ct
\tilde \to SV(V)(\tilde \leftarrow V^\infty(V).
\end{eqnarray*}

Let us consider the composition of $\gamma$ with $\ct$. Thus for any $0\leq l\leq n$
we get a map
\begin{eqnarray*}
\ct'_l\colon C^\infty(V\times \PP_+(V^*)^l,|\ome_V|\boxtimes
L^{\boxtimes l})\times C^\infty(V,|\ome_V|)\times
\ck(V)^{k+1}\times \RR^k_{\geq 0}\to \CC.
\end{eqnarray*}
Consider also the canonical $(l+1)$-linear map
\begin{eqnarray}\label{y1}
C^\infty(V,|\ome|)\times (C^\infty(\PP_+(V),L))^l\to
C^\infty(V\times \PP_+(V^*)^l,|\ome|\boxtimes L^{\boxtimes l}).
\end{eqnarray}
Composing $\ct'_l$ with the map (\ref{y1}) we get the map
\begin{eqnarray}\label{y2}
\ct_l''\colon C^\infty(V,|\ome|)\times (C^\infty(\PP_+(V^*),L))^l\times C^\infty(V,|\ome|)
\times \ck(V)^{k+1}\times \RR^k_{\geq 0}\to \CC
\end{eqnarray}
which is uniquely characterized by the following property:
\newline
for any $\mu,\nu\in C^\infty(V,|\ome_V|)$, any $B^1,\dots,B^l\in
\ck(V)$ being strictly convex compact sets with smooth boundaries,
and any $K,A^1,\dots,A^k\in \ck(V)$ one has
\begin{eqnarray*}
\ct_l''(\nu;h_{B^1},\dots,h_{B^l};\mu;K;A^1,\dots,A^k;\lam_1,\dots,\lam_k)=\\
\int_{x\in V} \frac{\pt^l}{\pt\mu_1\dots \pt\mu_l}\big|_0\nu\left(
(K\cap
(x-\sum_{i=1}^k\lam_iA^i))+\sum_{j=1}^l\mu_jB^j\right)d\mu(x).
\end{eqnarray*}

Using the L. Schwartz kernel theorem it is easy to see that in
order to prove Lemma \ref{l8.05} it is enough to show that for any
$0\leq l\leq n$ the map $\ct_l''$ has the following properties:

(1) $\ct_l''$ is $C^\infty$-smooth for fixed $K,A^1,\dots,A^k$;

(2) for any $R>0$ there exists a continous semi-norm $||\cdot||$
on $C^\infty(V,|\ome_V|)$ such that for any $K\in \ck(V)$
contained in the origin symmetric Euclidean ball of radius $R$,
and any strictly convex compact sets with smooth boundaries
$B_1,\dots,B_l\in\ck(V)$ one has an esimate
\begin{eqnarray}\label{y3}
\big| \frac{\pt^k}{\pt\lam_1\dots \pt\lam_k}\big|_0
\ct_l''(\nu;B_1,\dots,B_l;\mu;K,A_1,\dots,A_k;\lam_1,\dots,\lam_k)\big|\leq\\
||\nu||\cdot||\mu|| \cdot \prod_{i=1}^k||h_{A_i}||_{C^2(\PP_+(V^*))}\cdot
\prod_{j=1}^l||h_{B_j}||_{C^2(\PP_+(V^*))}.
\end{eqnarray}

In order to prove the last inequality let us observe that for fixed $\lam_1,\dots,\lam_k\geq 0$
\begin{eqnarray}\label{y4}
\int_{x\in V}\nu\left( K\cap(x-\sum_{i=1}^k\lam_iA_i)+\sum_{j=1}^l\mu_jB_j\right)d\mu(x)=
(\mu\boxtimes \nu)(\Delta (K)+(\sum_{i=1}^k\lam_iA_i,\sum_{j=1}^l\mu_jB_j))
\end{eqnarray}
where $\Delta\colon V\inj V\times V$ is the diagonal imbedding.
Now the inequality follows from Corollary \ref{cor-mink}.
Corollary \ref{cor-mink} implies also the smoothness of $\ct''_l$.
\qed


\begin{lemma}\label{l8.1}
Let $X\subset X'\subset V\times V$ be as at the beginning of this
section. Assume that $p_1(X)$ is convex. Let $A\in \ck(V)$. Let
\begin{eqnarray*}
\psi=(\Theta\circ\Xi_l)(\nu\otimes h_{B^1}\otimes \dots\otimes
h_{B^l})
\end{eqnarray*}
where $B^1,\dots,B^l$ are strictly convex compact sets with smooth
boundaries and containing the origin in their interiors. Let
$\mu,\nu$ be smooth densities on $V$. Then there exists $\eps>0$
depending on $X$ and $X'$ only such that the function $f\colon
[0,\eps]^l\to \CC$ defined by
$$f(\mu_1,\dots,\mu_l):= (\mu\boxtimes \nu)(X+(A\times
(\sum_{i=1}^l\mu_iB^i)))$$ is $C^\infty$-smooth. Moreover
\begin{eqnarray}\label{c1.1}
\frac{\pt^l}{\pt\mu_1\dots\pt\mu_l}\big |_0 f(\mu_1,\dots,\mu_l)=
\int_{x\in V} \left((\tilde p_1\circ \tilde
p_2^{-1})_*\psi\right)(p_1(X)\cap (x-A))d\mu(x).
\end{eqnarray}
\end{lemma}
{\bf Proof.} Let us choose $\eps>0$ such that for any $x\in V$ the
map $$N((X+A\times\{0\})\cap(\{x\}\times V))\times (0,1]\to V$$
given by $(p,n,t)\mapsto p+t\sum_{i=1}^l\mu_i\nabla h_{B^i}(n)$ is a
homeomorphism of $N((X+A\times\{0\})\cap(\{x\}\times
V))\times(0,1]$ onto its image
$$\left(p_2((X+A\times\{0\})\cap(\{x\}\times
V))+\sum_{i=1}^l\mu_iB^i\right)\backslash
\left(p_2((X+A\times\{0\})\cap(\{x\}\times V))\right)$$ for $0<
\mu_1,\dots,\mu_l\leq \eps$. Such an $\eps$ exists due to Lemmas
\ref{l7}, \ref{l2}, and Proposition \ref{l3}.

Let us denote
$$f(\mu_1,\dots,\mu_l):=
(\mu\boxtimes\nu)(X+(A\times(\sum_{i=1}^l\mu_iB^i))).$$ We have
\begin{eqnarray*}
f(\mu_1,\dots,\mu_l)=\\
\int_{x\in V}d\mu(x) \nu\left( ((X+A\times\{0\})\cap(\{x\}\times
V))+\{x\}\times (\sum_{i=1}^l\mu_iB^i)\right)=\\
\int_{x\in V}d\mu(x) \int_{N((X+A\times\{0\})\cap(\{x\}\times
V))\times[0,1]}
\left(\xi(\sum_{i=1}^l\mu_ih_{B^i},\cdot)\right)^*\nu=\\
\int_{x\in V}d\mu(x) \int_{N((X+A\times\{0\})\cap(\{x\}\times
V))\times[0,1]}\eta_{\mu_1,\dots,\mu_l}
\end{eqnarray*}
where
$\eta_{\mu_1,\dots,\mu_l}:=\left(\xi(\sum_{i=1}^l\mu_ih_{B^i},\cdot)\right)^*\nu$.
Consider the natural projection
$$q\colon V\times\PP_+(V^*)\times [0,1]\to V\times \PP_+(V^*).$$
Set
$$\tilde\eta_{\mu_1,\dots,\mu_l}:=q_*\eta_{\mu_1,\dots,\mu_l}\in
C^\infty(V\times\PP_+(V^*),\Ome^{n-1}\otimes p^*o)$$ (here $q_*$
denotes integration along the fibers). Using Lemma \ref{lift lemma} and the relation (\ref{intersection}), we compute
\begin{eqnarray}
f(\mu_1,\dots,\mu_l)=\\
\int_{x\in V}d\mu(x)\int_{N^*((X+(A\times\{0\}))\cap(\{x\}\times
V))}\tilde\eta_{\mu_1,\dots,\mu_l}=\\
\int_{x\in V}d\mu(x)\int_{(\overline{\tilde p_1 \circ \tilde p_2^{-1}})_* N^*(p_1(X)\cap(x-A))}\tilde\eta_{\mu_1,\dots,\mu_l}=\\
 \label{c1.2}\int_{x\in
V}d\mu(x)\int_{N^*(p_1(X)\cap(x-A))} \overline{(\tilde p_1\circ \tilde
p_2^{-1})}^* \tilde\eta_{\mu_1,\dots,\mu_l}.
\end{eqnarray}
Here $\overline{(\tilde p_1\circ \tilde p_2^{-1})}$ denotes the natural lift of the diffeomorphism
$\tilde p_1\circ \tilde p_2^{-1}$ to $V\times \mathbb P_+(V^*)$.
It is easy to see that the map $\RR^l\to
C^\infty(V\times\PP_+(V^*),\Ome^{n-1}\otimes p^*o)$ given by
$(\mu_1,\dots,\mu_l)\mapsto \tilde\eta_{\mu_1,\dots,\mu_l}$ is
$C^\infty$-smooth. This and Lemma \ref{l8.05}(1) imply the first
statement of the lemma.

Let us prove the second statement. Observe that for any compact
semi-convex (= positive reach) subset $Y\subset V$
$$\psi(Y)=\frac{\pt^l}{\pt\mu_1\dots\pt\mu_l}\big |_0
\int_{N(Y)} \tilde\eta_{\mu_1,\dots,\mu_l}.$$ Hence
$$\frac{\pt^l}{\pt\mu_1\dots\pt\mu_l} \big|_0
f(\mu_1,\dots,\mu_l)=\int_{x\in V}d\mu(x)\left((\tilde p_1\circ
\tilde p_2^{-1})_*\psi\right)(p_1(X)\cap (x-A)).$$
 \qed

\begin{lemma}\label{l9}
Let $\psi\in V^\infty(V)$ be a smooth valuation of the form
\begin{eqnarray}\label{c2}
\psi=(\Theta\circ \Xi_l)(\sum_{N=1}^\infty\nu_N\otimes
h_{B^1_N}\otimes \dots\otimes h_{B^l_N})
\end{eqnarray}
with $\{\nu_N\}\subset C^\infty(V,|\ome_V|)$ being smooth
densities, and $B_N^i\in \ck(V)$ being strictly convex compact
sets with smooth boundaries, containing the origin in the
interior, and such that for any compact subset $T\subset V$ and
any $L\in \NN$
\begin{eqnarray}\label{c3}
\sum_{N=1}^\infty
||\nu_N||_{C^L(T)}\prod_{i=1}^l||h_{B_N^i}||_{C^L(S^{n-1})}<\infty.
\end{eqnarray}

Let $A\in \ck(V)$. Let $\mu\in C^\infty(V,|\ome_V|)$. Let
$X\subset X'\subset V\times V$ be as in the beginning of this
section. Assume that $p_1(X)$ is convex.

Then the series
\begin{eqnarray}\label{c4}
\sum_{N=1}^\infty\frac{\pt^l}{\pt\mu_1\dots \pt\mu_l}\big |_0
(\mu\boxtimes \nu_N)(X+A\times (\sum_{i=1}^l\mu_iB_N^i))
\end{eqnarray}
converges absolutely and its sum is equal to
\begin{eqnarray}\label{c5}
\int_{x\in V} \left( (\tilde p_1\circ \tilde p_2^{-1})_*
\psi\right)(p_1(X)\cap(x-A))d\mu(x)
\end{eqnarray}
where as previously $\tilde p_1,\tilde p_2\colon X'\to V$ are the
restrictions of the projections $p_1,p_2$ to $X'$.
\end{lemma}
{\bf Proof.} If the sum in (\ref{c2}) is finite then the statement
follows immediately from Lemma \ref{l8.1}. Next let us observe
that the expression (\ref{c5}) is continuous with respect to
$\psi\in V^\infty(V)$. Hence it is enough to check that the series
(\ref{c4}) converges absolutely.

Let us denote
$$\psi_N:=(\Theta\circ \Xi_l)
(\nu_N\otimes h_{B^1_N}\otimes \dots\otimes h_{B^l_N}).$$ By Lemma
\ref{l8.1} we have
\begin{eqnarray}
\sum_{N=1}^\infty\frac{\pt^l}{\pt\mu_1\dots \pt\mu_l}\big |_0
(\mu\boxtimes \nu_N)(X+A\times (\sum_{i=1}^l\mu_iB_N^i))=\\
\sum_{N=1}^\infty \int_{x\in V} \left( (\tilde p_1\circ \tilde
p_2^{-1})_* \psi_N\right)(p_1(X)\cap(x-A))d\mu(x).
\end{eqnarray}
It follows from the assumption (\ref{c3}) that the series
$\sum_{N=1}^\infty \psi_N$ converges absolutely in $V^\infty (V)$.
Hence the series $\sum_{N=1}^\infty(\tilde p_1\circ \tilde
p_2^{-1})_* \psi_N$ converges absolutely in $V^\infty (p_1(X'))$.

\qed



\begin{proposition}\label{l10}

Let $X\subset X'\subset V\times V$ be as at the beginning of this
section. Assume moreover that $p_1(X)$ and $p_2(X)$ are convex.
Fix $\delta >0$. Then there exists $\eps_0>0$ depending on
$\delta$, the $C^2$-norm of the map $(p_2\circ
p_1^{-1})|_{p_1(X')}$ and its inverse, and the distance from $X$
to $\pt X'$ only such that the following properties are satisfied.
Let $\ca=(A^1,\dots,A^k)$ and $\cb=(B^1,\dots,B^l)$ be any $k$-
and $l$-tuples respectively of strictly convex compact subsets in
$V$ with smooth boundaries, principal curvatures between $\delta$
and $1/\delta$, and containing the origin in their interiors. Let
 $\mu, \, \nu$ be any smooth densities on $V$.

(1) Then the function
$$f_{\ca,\cb}(\lam_1,\dots,\lam_k;
\mu_1,\dots,\mu_l):=(\mu\boxtimes
\nu)(X+(\sum_{i=1}^k\lam_iA^i,\sum_{j=1}^l \mu_jB^j))$$ is
$C^\infty$-smooth for $(\lam_1,\dots,\lam_k; \mu_1,\dots,\mu_l)\in
[0,\eps]^{k+l}$ for any $0<\eps<\eps_0$ and such that
$p_1(X)+\eps\sum_iA^i\subset p_1(X')$ and $p_2(X)+\eps\sum_j
B^j\subset p_2(X').$

(2) There exist continuous semi-norm $||\cdot ||$ on
$C^\infty(V,|\ome_V|)$ and a positive integer $L\in \NN$ depending
on $X,X',k,l$ only ( and independent of $\ca,\cb,\mu,\nu$) such
that
\begin{eqnarray*}
\big|
\frac{\pt^{k+l}}{\pt\lam_1\dots\pt\lam_k\pt\mu_1\dots\pt\mu_l}\big|_0
f_{\ca,\cb}(\lam_1,\dots,\lam_k; \mu_1,\dots,\mu_l)\big|\leq\\
C||\mu||\cdot ||\nu||\cdot
\prod_{i=1}^k||h_{A^i}||_{C^L(S^{n-1})}\cdot\prod_{j=1}^l||h_{B^j}||_{C^L(S^{n-1})}.
\end{eqnarray*}

\end{proposition}
{\bf Proof.} Similarly to the proof of Lemma \ref{l8.1}, let us
denote
\def\ebm{\eta_{\cb,\mu_1,\dots,\mu_l}}
\def\tebm{\tilde\eta_{\cb,\mu_1,\dots,\mu_l}}
\begin{eqnarray}
\eta_{\cb,\mu_1,\dots,\mu_l}:=(\xi(\sum_{i=1}^l\mu_i h_{B^i},
\cdot))^*\nu\in C^\infty(V\times \PP_+(V^*)\times
[0,1],\Ome^n\otimes p^*o)\\
\tilde
\eta_{\cb,\mu_1,\dots,\mu_l}:=q_*\eta_{\cb,\mu_1,\dots,\mu_l}\in
C^\infty(V\times \PP_+(V^*),\Ome^{n-1}\otimes p^*o),
\end{eqnarray}
where as previously $q\colon V\times\PP_+(V^*)\times [0,1]\to
V\times \PP_+(V^*)$ is the projection, and $q_*$ denotes the
integration along the fibers.

\def\fab{f_{\ca,\cb}}
\def\fabm{f_{\ca,\cb}(\lam_1,\dots,\lam_k;\mu_1,\dots,\mu_l)}
First let us prove part (1) of the proposition. Exactly as in
(\ref{c1.2}) we have
\begin{eqnarray*}
\fabm=\int_{x\in V} d\mu(x)\int_{N(p_1(X)\cap
(x-\sum_i\lam_iA^i))}\overline{(\tilde p_2\circ\tilde p_1^{-1})}^*\tebm.
\end{eqnarray*}
Lemma \ref{l4} implies that $\tebm$ depends smoothly on
$(\mu_1,\dots,\mu_l)\in [0,\eps]^l$ and on $h_{B^1},\dots,h_{B^l}
\in C^\infty(S^{n-1})$.
\def\zbm{\zeta_{\cb,\mu_1,\dots,\mu_l}}
Hence
\begin{eqnarray*}
\zbm:=\overline{(\tilde p_2\circ\tilde p_1^{-1})}^*\tebm
\end{eqnarray*}
also depends smoothly on $(\mu_1,\dots,\mu_l)\in [0,\eps]^l$ and
on $h_{B^1},\dots,h_{B^l}\in C^\infty(S^{n-1})$. Thus
\begin{eqnarray*}
\fabm=\int_{x\in V}d\mu(x)\int_{N(p_1(X)\cap
(x-\sum_i\lam_iA^i))}\zbm.
\end{eqnarray*}
Then by Lemma \ref{l8.05}(1) the function $\fab$ is
$C^\infty$-smooth in $(\lam_1,\dots,\lam_k;\mu_1,\dots,\mu_l)\in
[0,\eps]^{k+l}$. This proves part (1) of the proposition.

Let us prove part (2).
\def\sib{\sigma_\cb}
Let us denote by
\begin{eqnarray*}
\sib:=\frac{\pt^l}{\pt\mu_1\dots \pt\mu_l}\big|_0\zbm.
\end{eqnarray*}
Then $\sib\in C^\infty(V\times\PP_+(V^*),\Ome^{n-1}\otimes p^*o)$.
Moreover $\sib$ depends continuously and linearly on each
$h_{B^i}\in C^\infty(S^{n-1})$. Then for any $M\in \NN$ there
exist a compact subset $T'\subset V$, $L\in \NN$, and a constant
$C$ such that
\begin{eqnarray}\label{t1}
||\sigma_\cb||_{C^M(p_1(X'))}\leq C||\nu||_{C^L(T')}\cdot
\prod_{j=1}^l||h_{B^j}||_{C^L(S^{n-1})}.
\end{eqnarray}
We have
\begin{eqnarray*}
\frac{\pt^l}{\pt\mu_1\dots\pt\mu_l}\big|_0\fabm= \int_{x\in
V}d\mu(x) \int_{N(p_1(X)\cap (x-\sum_i\lam_iA^i))}\sib
\end{eqnarray*}
(note that the differentiation under the the integral is possible
due to Lemma \ref{l8.05}(1)).

Hence by Lemma \ref{l8.05}(2) there exist continuous semi-norms
$||\cdot||,\, ||\cdot ||'$ on $C^\infty(p_1(X')\times
\PP_+(V^*),\Ome^{n-1}\otimes p^*o)$ and $C^\infty(V,|\ome_V|)$
respectively and $L'\in \NN$ such that
\begin{eqnarray}
\big|\frac{\pt^{k+l}}{\pt\lam_1\dots\lam_k\pt\mu_1\dots\pt\mu_l}\big
|_0\fabm\big|\leq \\
||\sib||\cdot ||\mu||'\cdot
\prod_{i=1}^k||h_{A^i}||_{C^{L'}(S^{n-1})} \overset{\mbox{by }
(\ref{t1})}{\leq} \\
\label{t2} C||\nu||_{C^L(T')}\cdot ||\mu||'\cdot \prod _{i=1}^k
||h_{A^i}||_{C^{L'}(S^{n-1})}\cdot
\prod_{j=1}^l||h_{B^j}||_{C^L(S^{n-1})}.
\end{eqnarray}
Note that the semi-norms $||\cdot||,||\cdot||'$, and the constant
$C$ in (\ref{t2}) are independent of $\ca,\cb$. This proves part
(2) of the proposition. \qed



\begin{lemma}\label{l11}
Let $Y$ be a smooth $n$-dimensional manifold. Let $f_1,f_2\colon
Y\to \RR^n$ be two smooth maps which map $Y$ diffeomorphically
onto open subsets $f_1(Y),f_2(Y)\subset\RR^n$. Let $\phi\in
V^\infty(Y)$. Assume that $\phi(K)=0$ for any compact domain
$K\subset Y$ with smooth boundary  such that both $f_1(K)$ and
$f_2(K)$ are convex.

Then $\phi\equiv 0$.
\end{lemma}
{\bf Proof.} Suppose that $\phi\not\equiv 0$. Let $i$ be the
integer such that $\phi\in W_i\backslash W_{i+1}$.

Let $K\subset Y$ be a compact domain with smooth boundary such
that $f_1(K)\subset \RR^n$ is strictly convex. Let us show that
for any $y_0\in Y$ the set $f_2(f_1^{-1}(f_1(y_0)+\eps f_1(K)))$
is convex for $0< \eps \ll 1$. For simplicity of the notation and
without loss of generality we may assume that $Y=f_1(Y)$ and $f_1$
is the identity imbedding. Also we may and will assume that
$f_2(y)=y+O(|y-y_0|^2)$ as $y\to y_0$, and $y_0=0$. Then
$$\frac{f_2(\eps y)}{\eps}=y+O(\eps)\mbox{ as } \eps\to 0.$$
But since $K$ is strictly convex it is clear that from the last
formula that $\frac{f_2(\eps K)}{\eps}$ is also convex for
$0<\eps\ll 1$.

Now let us deduce the lemma. Fix $K\subset Y$ a compact domain
with smooth boundary such that $f_1(K)$ is strictly convex. For
$0<\eps \ll 1$ the set $f_2(f_1^{-1}(f_1(y_0)+\eps f_1(K)))$ is
convex. Hence the assumption of the lemma imply that
\begin{eqnarray*}
\lim_{\eps\to +0} \frac{1}{\eps^i} \phi(f_1^{-1}(f_1(y_0)+\eps
f_1(K)))=0.
\end{eqnarray*}
Hence $\phi\in W_{i+1}$. This is a contradiction. \qed


\section{Construction of the product.}\label{construction}
\setcounter{subsection}{1}
\setcounter{theorem}{0}

Let $X$ be a smooth
$n$-dimensional manifold. Let $\phi,\psi\in V^\infty(X)$. We are
going to present a construction of a product $\phi\cdot\psi$ which
{\itshape a priori} will depend on some choices, and then we will
show that the product is in fact independent of these choices.

Let $U\subset X$ be an open subset. Let $f\colon U\tilde \to
\RR^n$ be a diffeomorphism. It was shown in \cite{part1}, Corollary 3.1.7, that $f_*\phi$ can
be written (non-uniquely) as
\begin{eqnarray}\label{a1}
f_*\phi=\phi_0+\dots+\phi_n
\end{eqnarray}
with $\phi_j\in W_{n-j}(\RR^n)$, and there exist sequences
\begin{eqnarray}\label{a2}
\{\mu_N^j\}_{N=1}^\infty\subset C^\infty(\RR^n,|\ome_{\RR^n}|),\,
0\leq j\leq n;\\
\label{a2'}\{A_N^{ij}\}_{N=1}^\infty\subset \ck(\RR^n),\, 1\leq
i\leq j\leq n
\end{eqnarray}
being strictly convex compact domains with smooth boundaries and
containing 0 in their interiors such that

($1_{f_*\phi}$) for any compact subset $T\subset \RR^n$, any
$L\in\NN$, and any $0\leq j\leq n$ one has
\begin{eqnarray}\label{a3}
\sum_{N=1}^\infty ||\mu_N^j||_{C^L(T)}\prod _{i=1}^j
||h_{A_N^{ij}}||_{C^L(S^{n-1})}< \infty;
\end{eqnarray}

($2_{f_*\phi}$) for any set $S\in \ck(\RR^n)$ and any $0\leq j\leq
n$ one has
\begin{eqnarray}\label{a4}
\phi_j(S)=\sum_{N=1}^\infty \frac{\pt^j}{\pt \lam_1\dots
\pt\lam_j}\big |_0 \mu_N^j(S+\sum_{i=1}^j\lam_iA_N^{ij}).
\end{eqnarray}
The last expression is well defined since by Corollary \ref{cor-mink} the function
$\mu_N^j(S+\sum_{i=1}^j\lam_iA_N^{ij})$ is $C^\infty$-smooth in
$\lam_1,\dots,\lam_j\geq 0$ and the series (\ref{a4}) converges
absolutely.

Similarly one can write

\begin{eqnarray}\label{a5}
f_*\psi=\psi_0+\dots+\psi_n
\end{eqnarray}
with $\psi_j\in W_{n-j}(\RR^n)$, and there exist sequences
\begin{eqnarray}\label{a6}
\{\nu_N^j\}_{N=1}^\infty\subset C^\infty(\RR^n,|\ome_{\RR^n}|),\,
0\leq j\leq n;\\
\label{a6'}\{B_N^{ij}\}_{N=1}^\infty\subset \ck(\RR^n),\, 1\leq
i\leq j\leq n
\end{eqnarray}
with (\ref{a6'}) being strictly convex compact domains with smooth
boundaries containing 0 at the interior such that

($1_{f_*\psi}$) for any compact subset $T\subset \RR^n$, any
$L\in\NN$, and any $0\leq j\leq n$ one has
\begin{eqnarray}\label{a7}
\sum_{N=1}^\infty ||\nu_N^j||_{C^L(T)}\prod _{i=1}^j
||h_{B_N^{ij}}||_{C^L(S^{n-1})}< \infty;
\end{eqnarray}

($2_{f_*\psi}$) for any set $S\in \ck(\RR^n)$ and any $0\leq j\leq
n$ one has
\begin{eqnarray}\label{a8}
\psi_j(S)=\sum_{N=1}^\infty \frac{\pt^j}{\pt \mu_1\dots
\pt\mu_j}\big |_0 \nu_N^j(S+\sum_{i=1}^j\mu_i B_N^{ij}).
\end{eqnarray}
As previously the function $\nu_N^j(S+\sum_{i=1}^j\mu_i B_N^{ij})$
is $C^\infty$-smooth in $\mu_1,\dots,\mu_j\geq 0$ and the series
(\ref{a8}) converges absolutely.

\def\pjj{\frac{\pt^{j+j'}}{\pt\lam_1\dots\pt\lam_j\pt\mu_1\dots\pt\mu_{j'}}\big
|_0} In \cite{part1} we have defined the product $f_*\phi\cdot
f_*\psi$ as a valuation defined on {\itshape convex} compact
subsets of $\RR^n$ only by the following formula: for any $K\in
\ck(\RR^n)$
\begin{eqnarray}\label{a10}
(f_*\phi\cdot f_*\psi)(K)=\\
\label{a10.1}\sum_{j,j'=0}^n\sum_{N,N'=1}^\infty \pjj
(\mu_N^j\boxtimes \nu_{N'}^{j'})\left(\Delta
(K)+(\sum_{i=1}^j\lam_i A_N^{ij}\times
\sum_{i'=1}^{j'}\mu_{i'}B_{N'}^{i'j'})\right)
\end{eqnarray}
where $\Delta \colon \RR^n\hookrightarrow \RR^n\times \RR^n$ is
the diagonal imbedding, the function $$(\mu_N^j\boxtimes
\nu_{N'}^{j'})\left(\Delta (K)+(\sum_{i=1}^j\lam_i A_N^{ij}\times
\sum_{i'=1}^{j'}\mu_{i'}B_{N'}^{i'j'})\right)$$ is
$C^\infty$-smooth in $\lam_1,\dots,\lam_j,\mu_1,\dots,\mu_{j'}\geq
0$, the series (\ref{a10.1}) converges absolutely and defines a
valuation on $\ck(\RR^n)$ from the space $SV(\RR^n)$. By
Proposition \ref{isomorphism} it defines a smooth valuation on
$\cp(\RR^n)$. Hence we get a smooth valuation on $U\subset X$.
This valuation will be denoted later on by $\phi|_U\cdot \psi|_U$.
However this construction depends {\itshape a priori} on a choice
of a diffeomorphism $f$ and the choices (\ref{a1})-(\ref{a2'}),
(\ref{a5})-(\ref{a6'}). It was however shown in \cite{part1} that
once $f$ is fixed, the other choices (\ref{a1})-(\ref{a2'}),
(\ref{a5})-(\ref{a6'}) do not influence the definition of
$f_*\phi\cdot f_*\psi$. So let us denote temporarily the valuation
we have constructed on $U$ by $\phi|_U\circ_f \psi|_U$. In order
to check that the product is well defined it remains to show that
if $\tilde U\subset X$ is another open subset and $\tilde f\colon
\tilde U \tilde \to \RR^n$ is a diffeomorphism then
\begin{eqnarray}\label{a13}
(\phi|_U\circ_f\psi|_U)|_{U\cap\tilde U}=(\phi|_{\tilde
U}\circ_{\tilde f}\psi|_{\tilde U})|_{U\cap \tilde U}.
\end{eqnarray}

By Lemma \ref{l11} it is enough to show that for an arbitrary
compact domain with smooth boundary $K\subset U\cap \tilde U$ such
that $f(K)$ and $\tilde f(K)$ are convex in $\RR^n$, one has
\begin{eqnarray}\label{a14}
(f_*\phi\cdot f_*\psi)(f(K))=(\tilde f_*\phi\cdot \tilde f
_*\psi)(\tilde f(K)).
\end{eqnarray}
\def\tph{\tilde f_*\phi}
\def\tps{\tilde f_*\psi}
Let us also fix another compact domain with smooth boundary
$K'\subset U\cap \tilde U$ such that $K$ is contained in the
interior of $K'$. For valuations $\tph$ and $\tps$ we can find
presentations similar to (\ref{a1}), (\ref{a4}), (\ref{a5}),
(\ref{a8}). Namely
\begin{eqnarray}
\label{a15} \tph=\tilde \phi_0+\dots +\tilde \phi_n,\\
\label{a16} \tps=\tilde \psi_0+\dots +\tilde \psi_n,
\end{eqnarray}
with $\tilde\phi_{j},\tilde \psi_j \in W_{n-j}(\RR^n)$ and there
exist sequences
\begin{eqnarray}\label{a17}
\{\tilde\mu_N^j\}_{N=1}^\infty,\,\{\tilde\nu_N^j\}_{N=1}^\infty
\subset C^\infty(\RR^n,|\ome_{\RR^n}|),\,
0\leq j\leq n;\\
\label{a18}\{\tilde A_N^{ij}\}_{N=1}^\infty,\,\{\tilde
B_N^{ij}\}_{N=1}^\infty \subset \ck(\RR^n),\, 1\leq i\leq j\leq n
\end{eqnarray}
with (\ref{a18}) being strictly convex compact domains with smooth
boundaries, containing the origin at the interior, and such that

(1) for any compact subset $T\subset \RR^n$, any $L\in\NN$, and
any $0\leq j\leq n$ one has
\begin{eqnarray}\label{a19}
\sum_{N=1}^\infty ||\tilde \mu_N^j||_{C^L(T)}\prod _{i=1}^j
||h_{\tilde A_N^{ij}}||_{C^L(S^{n-1})}< \infty,\\
\label{a20} \sum_{N=1}^\infty ||\tilde \nu_N^j||_{C^L(T)}\prod
_{i=1}^j ||h_{\tilde B_N^{ij}}||_{C^L(S^{n-1})}< \infty;
\end{eqnarray}

(2) for any set $S\in \ck(\RR^n)$ and any $0\leq j\leq n$ one has
\begin{eqnarray}\label{a21}
\tilde \phi_j(S)=\sum_{N=1}^\infty \frac{\pt^j}{\pt \lam_1\dots
\pt\lam_j}\big |_0 \tilde \mu_N^j(S+\sum_{i=1}^j\lam_i\tilde A_N^{ij}),\\
\label{a22}\tilde \psi_j(S)=\sum_{N=1}^\infty \frac{\pt^j}{\pt
\mu_1\dots \pt\mu_j}\big |_0 \tilde
\nu_N^j(S+\sum_{i=1}^j\mu_i\tilde B_N^{ij}) .
\end{eqnarray}

We have to check that
\begin{eqnarray}\label{a23}
\sum_{j,j'=0}^n\sum_{N,N'=1}^\infty \pjj (\mu_N^j\boxtimes
\nu_{N'}^{j'})\left((f\times f) (K)+(\sum_{i=1}^j\lam_i
A_N^{ij}\times \sum_{i'=1}^{j'}\mu_{i'}B_{N'}^{i'j'})\right)=\\
\label{a23.5} \sum_{j,j'=0}^n\sum_{N,N'=1}^\infty \pjj
(\tilde\mu_N^j\boxtimes \tilde\nu_{N'}^{j'})\left((\tilde f\times
\tilde f) (K)+(\sum_{i=1}^j\lam_i \tilde A_N^{ij}\times
\sum_{i'=1}^{j'}\mu_{i'}\tilde B_{N'}^{i'j'})\right).
\end{eqnarray}
We will prove the following lemma.
\begin{lemma}\label{lemma1}
The expression (\ref{a23}) is equal to
\begin{eqnarray}\label{a24}
\sum_{j,j'=0}^n\sum_{N,N'=1}^\infty \pjj (\mu_N^j\boxtimes \tilde
\nu_{N'}^{j'})\left((f\times \tilde f) (K)+(\sum_{i=1}^j\lam_i
A_N^{ij}\times \sum_{i'=1}^{j'}\mu_{i'}\tilde
B_{N'}^{i'j'})\right).
\end{eqnarray}
In the last expression (\ref{a24}) the function $(\mu_N^j\boxtimes
\tilde \nu_{N'}^{j'})\left((f\times \tilde f)
(K)+(\sum_{i=1}^j\lam_i A_N^{ij}\times
\sum_{i'=1}^{j'}\mu_{i'}\tilde B_{N'}^{i'j'})\right)$ is
$C^\infty$-smooth in $0\leq
\lam_1,\dots,\lam_j,\mu_1,\dots,\mu_{j'}< \eps_N$ for some
$\eps_N>0$ depending on $K,A_N^{ij},\tilde B_{N'}^{i'j'}$, and the
series converges absolutely.
\end{lemma}

Let us show first that Lemma \ref{lemma1} implies the equality
(\ref{a14}) and hence implies that the product of valuations is
well defined. We can apply Lemma \ref{lemma1} once again in a
symmetric way in order to show that the expression (\ref{a24}) is
equal to the expression (\ref{a23.5}). Thus the equality
(\ref{a23})=(\ref{a23.5}) will be proved.

{\bf Proof of Lemma \ref{lemma1}.} The differentiability and
absolute convergence in (\ref{a24})  follow immediately from
Proposition \ref{l10}.


Let us show that (\ref{a23})=(\ref{a24}). Let us fix $j,N$. It is
enough to show that
\begin{eqnarray*}
\sum_{j'=0}^n\sum_{N'=1}^\infty\frac{\pt^{j+j'}}{\pt\lam_1\dots\pt\lam_j\pt\mu_1\dots\pt\mu_{j'}}\big|_0
(\mu_N^j\boxtimes \nu_{N'}^{j'})((f\times
f)(K)+(\sum_{i=1}^j\lam_iA^{ij}_N\times \sum_{i'=1}^{j'}\mu_{i'}B^{i'j'}_{N'}))=\\
\sum_{j'=0}^n\sum_{N'=1}^\infty\frac{\pt^{j+j'}}{\pt\lam_1\dots\pt\lam_j\pt\mu_1\dots\pt\mu_{j'}}\big|_0
(\mu_N^j\boxtimes \tilde\nu_{N'}^{j'})((f\times \tilde
f)(K)+(\sum_{i=1}^j\lam_iA^{ij}_N\times
\sum_{i'=1}^{j'}\mu_{i'}\tilde B^{i'j'}_{N'})).
\end{eqnarray*}
Lemma \ref{l8.05} implies that it is enough to show that for fixed
$\lam_1,\dots,\lam_j>0$ one has the equality
\begin{eqnarray*}
\sum_{j'=0}^n\sum_{N'=1}^\infty\frac{\pt^{j'}}{\pt\mu_1\dots\pt\mu_{j'}}\big|_0
(\mu_N^j\boxtimes \nu_{N'}^{j'})((f\times
f)(K)+(\sum_{i=1}^j\lam_iA^{ij}_N\times\sum_{i'=1}^{j'}\mu_{i'}B^{i'j'}_{N'}))=\\
\sum_{j'=0}^n\sum_{N'=1}^\infty\frac{\pt^{j'}}{\pt\mu_1\dots\pt\mu_{j'}}\big|_0
(\mu_N^j\boxtimes \tilde\nu_{N'}^{j'})((f\times \tilde
f)(K)+(\sum_{i=1}^j\lam_iA^{ij}_N\times\sum_{i'=1}^{j'}\mu_{i'}\tilde
B^{i'j'}_{N'})).
\end{eqnarray*}
Let us denote for brevity $A:=\sum_{i=1}^j\lam_iA^{ij}_N,\,
\mu:=\mu_N^j$. In this notation the last equality is rewritten
\begin{eqnarray*}
\sum_{j'=0}^n\sum_{N'=1}^\infty\frac{\pt^{j'}}{\pt\mu_1\dots\pt\mu_{j'}}\big|_0
(\mu\boxtimes \nu_{N'}^{j'})((f\times
f)(K)+(A\times\sum_{i'=1}^{j'}\mu_{i'}B^{i'j'}_{N'}))=\\
\sum_{j'=0}^n\sum_{N'=1}^\infty\frac{\pt^{j'}}{\pt\mu_1\dots\pt\mu_{j'}}\big|_0
(\mu\boxtimes \tilde\nu_{N'}^{j'})((f\times \tilde
f)(K)+(A\times\sum_{i'=1}^{j'}\mu_{i'}\tilde B^{i'j'}_{N'})).
\end{eqnarray*}
By Lemma \ref{l9} both sides of the last equality are equal to
$$\int_{x\in V} ((f\circ\tilde
f^{-1})_*\psi)(f(K)\cap(x-A))d\mu(x).$$ This finishes the proof of
Lemma \ref{lemma1}. Hence this also finishes the proof that the
product on smooth valuations is well defined. \qed

>From the construction of the product it is easy to deduce the
following result.
\begin{proposition}\label{q1}
Let $U\subset V\subset X$ be open subsets of $X$. Let $\phi,\psi\in V^\infty(V)$.
Then $$(\phi\cdot \psi)|_U=\phi|_U\cdot \psi|_U.$$
\end{proposition}

\begin{theorem}\label{q2}
The product
$$V^\infty(X)\times V^\infty(X)\to V^\infty(X)$$
is continuous, commutative, and associative. The Euler characteristic is the
unit in the algebra $V^\infty(X)$.
\end{theorem}

{\bf Proof.} Observe first that if $X$ is diffeomorphic to $\RR^n$
then this theorem was proved in \cite{part1}, Theorem 4.1.2
(combined with the description of $V^\infty(\RR^n)$ from
Proposition 2.4.10 in \cite{part2}). This and Proposition \ref{q1}
imply all the statements of the theorem except of continuity.

Let us prove continuity. Assume that $\phi_N\to \phi,\psi_N\to\psi$ in
$V^\infty(X)$. We have to show that $\phi_N\cdot\psi_N\to \phi\cdot\psi$ in
$V^\infty(X)$. Note that for any open subset $U\subset X$ diffeomeorphic to $\RR^n$
$$(\phi_N\cdot \psi_N)|_U\to (\phi\cdot \psi)|_U \mbox{ in } V^\infty(U)$$
by the affine case and Proposition \ref{q1}.

One can easily check the following property. Let $\{\xi_N\}\subset V^\infty(X),\,
\xi\in V^\infty(X)$. Let $\{U_\alp\}_\alp$ be an open covering of $X$. Assume that
for any $\alp$
$$\xi_N|_{U_\alp}\to \xi|_{U_\alp} \mbox{ in } V^\infty(U_\alp).$$
Then $\xi_N\to \xi$ in $V^\infty(X)$.
This implies the theorem. \qed

Recall now that by \cite{part2} the assignment to any open subset $U\subset X$
$$U\mapsto V^\infty(U)$$
is a sheaf on $X$ denoted by $\cv^\infty_X$. Proposition \ref{q1} and Theorem \ref{q2}
imply immediately the following corollary.
\begin{corollary}\label{q3}
$\cv^\infty_X$ is a sheaf of commutative associative algebras with unit
(where the unit is the Euler characteristic).
\end{corollary}




\noindent {\sc Semyon Alesker\\
{ \normalsize Department of Mathematics, Tel Aviv University,
Ramat Aviv}
 \\  { \normalsize 69978 Tel Aviv,
Israel }
\\ \tt semyon@post.tau.ac.il\\

\noindent \sc Joseph H.G. Fu\\
{\normalsize Department of Mathematics, University of Georgia}\\
{\normalsize Athens, Georgia 30602, USA}\\
\tt fu@math.uga.edu

\end{document}